\documentclass[12pt,A4paper]{article}

\usepackage[T1]{fontenc}
\usepackage[latin1]{inputenc}

\usepackage{amssymb}
\makeatletter
\@addtoreset{equation}{section}

\makeatother

\newtheorem{Th}{Theorem}[section]
\newtheorem{Theo}[Th]{Theorem}
\newtheorem{thm}[Th]{Theorem}
\newtheorem{Prop}[Th]{Proposition}
\newtheorem{prop}[Th]{Proposition}

\newtheorem{cor}[Th]{Corollary}

\newtheorem{lem}[Th]{Lemma}

\newtheorem{rem}[Th]{Remark}
\newtheorem{example}[Th]{Example}

\begin{document}

\title{Asymptotics for Bergman-Hodge kernels for high powers of 
complex line bundles}
\author{Robert Berman\\Department of Mathematics, Chalmers University
  of Technology,\\Eklandag.
86, SE-412 96 G\"oteborg \and 
Johannes Sj\"ostrand\\CMLS, Ecole Polytechnique, FR-91128 Palaiseau 
cedex,\\ UMR
    7640, CNRS
}
\date{} 
\maketitle
\newcommand{\iint}{\int\hskip -2mm\int}
\newcommand{\iiint}{\int\hskip -2mm\int\hskip -2mm\int}
\def\ekv#1#2{\begin{equation}\label{#1}#2\end{equation}}
\def\eekv#1#2#3{\begin{eqnarray}\label{#1}#2\\ #3\nonumber\end{eqnarray}}
\def\no#1{(\ref{#1})}
\def\cint{{1\over 2\pi i}\int}
\def\iint{\int\hskip -2mm\int}
\def\ciint{({1\over 2\pi i})^2\iint}
\def\iiint{\int\hskip -2mm\int\hskip -2mm\int}
\def\iiiint{\int\hskip -2mm\int\hskip -2mm\int \hskip -2mm\int}
\def\buildover#1#2{\buildrel#1\over#2}
\font \mittel=cmbx10 scaled \magstep1
\font \gross=cmbx10 scaled \magstep2
\font \klo=cmsl8
\font\liten=cmr10 at 8pt
\font\stor=cmr10 at 12pt
\font\Stor=cmbx10 at 14pt
\def\aby{arbitrary}
\def\ably{arbitrarily}
\def\an{analytic}
\def\asy{asymptotic}
\def\bdd{bounded}
\def\bdy{boundary}
\def\coef{coefficient}
\def\const{constant}
\def\Const{Constant}
\def\canform{canonical transformation}
\def\cont{continous}
\def\ctf{canonical transformation}
\def\diffeo{diffeomorphism}
\def\de{differential equation}
\def\dop{differential operator}
\def\ef{eigenfunction}
\def\ev{eigenvalue}
\def\e{equation}
\def\fu{function}
\def\fy{family}
\def\F{Fourier}
\def\fop{Fourier integral operator}
\def\fourior{Fourier integral operator}
\def\fouriors{Fourier integral operators }
\def\hol{holomorphic}
\def\indep{independent}
\def\lhs{left hand side}
\def\mfld{manifold}
\def\ml{microlocal}
\def\neigh{neighborhood}
\def\nondeg{non-degenerate}
\def\op{operator}
\def\og{orthogonal}
\def\pb{problem}
\def\Pb{Problem}
\def\pde{partial differential equation}
\def\pe{periodic}
\def\pert{perturbation}
\def\plsh{plurisubharmonic}
\def\plh{pluriharmonic}
\def\pol{polynomial}
\def\pro{proposition}
\def\Pro{Proposition}
\def\pop{pseudodifferential operator}
\def\Pop{Pseudodifferential operator}
\def\pseudor{pseudodifferential operator}
\def\res{resonance}
\def\rhs{right hand side}
\def\sa{selfadjoint}
\def\schr{Schr\"odinger operator}
\def\scl{semiclassical}
\def\sop{Schr\"odinger operator}
\def\st{strictly}
\def\stpsh{\st{} plurisubharmonic}
\def\suf{sufficient}
\def\sufly{sufficiently}
\def\tf{transformation}
\def\Th{Theorem}
\def\th{theorem}
\def\tf{transform}
\def\traj{trajectory}
\def\trans{^t\hskip -1pt }
\def\top{Toeplitz operator}
\def\uf{uniform}
\def\ufly{uniformly}
\def\vf{vector field}
\def\wrt{with respect to}
\def\Op{{\rm Op\,}}
\def\Re{{\rm Re\,}}
\def\Im{{\rm Im\,}}
\def\dbar{\overline{\partial }}
\def\rf{\rfloor}
\def\db{\overline{\partial }}
\def\on{orthonormal}
\def\C{{\bf C}}
\def\R{{\bf R}}
\def\Z{{\bf Z}}
\begin{abstract}
In this paper we obtain the full \asy{} expansion of the Bergman-Hodge kernel 
associated to a high power of a \hol{} line bundle with non-degenerate 
curvature. We also explore some relations with \asy{} \hol{} sections on 
symplectic \mfld{}s.\medskip\par\noindent
\centerline{\hskip -20mm \bf R\'esum\'e}\smallskip
\par Dans ce travail nous obtenons un d\'eveloppement asymptotique complet 
du noyau de Bergman-Hodge d'une puissance \'elev\'ee d'un fibr\'e en droites 
holomorphe \`a courbure non-d\'eg\'ener\'ee. Nous explorons aussi quelques 
relations avec des sections asymptotiquement holomorphes sur une 
vari\'et\'e symplectique.
\end{abstract}
\tableofcontents
\section{Introduction}\label{int}
Let $L$ be a Hermitian holomorphic line bundle over a compact complex
Hermitian manifold $X.$ Denote by $\Theta$ the curvature two-form
of the canonical connection $\nabla$ on $L.$ By the Hodge theorem,
the Dolbeault cohomology group $H^{0,q}(X,L)$ is isomorphic to the
space $\mathcal{H}^{0,q}(X,L)$ of harmonic $(0,q)-$forms with values
in $L,$ i.e the null space of the Hodge Laplacian $\Delta_{q}$.
Denote by $\Pi_{q}$ the corresponding Hodge projection, i.e. the
orthogonal projection from $L^{2}(X,L)$ onto $\mathcal{H}^{0,q}(X,L).$
We will assume that $\Theta$ is non-degenerate of constant signature
$(n_{-},n_{+})$, i.e. the number of negative eigenvalues of $\Theta$
is $n_{-}$ (the index of $\Theta).$ Then it is well-known, by the
theorems of Kodaira and Hörmander, that $\mathcal{H}^{0,q}(X,L^{k})$
is trivial when $q\neq n_{-},$ for a sufficiently high tensor power
$L^{k}$. (See also \cite{DePeSc}.) We will study the asymptotics with respect to $k$ of the
corresponding Hodge projections $\Pi_{q,k}$ in the non-trivial case
when $q=n_{-}.$ The case when $n_{-}=0,$ i.e. when $L$ is a positive
line bundle and $\Pi_{q,k}$ is the Bergman projection on the space
of holomorphic sections with values in $L^{k},$ has been studied
extensively before (compare the historical remarks below). 

Let $\pi_{1}$ and $\pi_{2}$ be the projections on the first and
the second factor of $X\times X.$ Denote by $K_{k}$ the Schwartz
kernel of $\Pi_{q,k}$ (the subscripts $k$ will be omitted in the
sequel) with respect to the volume form $\omega_{n}$ on $X$ induced
by the Hermitian metric on $X,$ so that $K$ is a section of 
${\cal L}(\pi_2^*(\Lambda^{0,q}(T^*X\otimes 
L^{k})),\pi_1^*(\Lambda^{0,q}(T^*X\otimes L^{k})))$.

Let $t,s$ be local unitary sections of $L$ over 
$\widetilde{X}$, $\widetilde{Y}$
respectively, where $\widetilde{X},\widetilde{Y}\subseteq X.$ Then
on $\widetilde{X}\times\widetilde{Y}$ we can write \[
K(x,y)=K_{t,s}(x,y;\frac{1}{k})t(x)^{k}s(y)^{*k},\]
where $K_{t,s}$ is a local section of 
$\mathcal{L}(\textrm{$\pi_{2}^{*}(\Lambda^{0,q}(T^{*}X)),\pi_{1}^{*}
(\Lambda^{0,q}(T^{*}X)))$ 
}$so that for 
$x\in\widetilde{X}$, $u\in C_{0}^{\infty}(\widetilde{Y};
\Lambda^{0,q}(T^{*}X\otimes L^{k})),$
\[
u(x)=t(x)^{k}\int_{X}K_{t,s}(x,y;\frac{1}{k})\left\langle 
u(y),s(y)^{*k}\right\rangle \omega_{n}(dy),\]
We say that a kernel \[
R(x,y)=R_{t,s}(x,y;\frac{1}{k})t(x)^{k}s(y)^{*k},\]
 is \emph{negligible} if \[
\partial_{x}^{\alpha}\partial_{y}^{\beta}R_{t,s}(x,y;\frac{1}{k})=\mathcal{O}_{\alpha,\beta,N}(k^{-N}),\]
locally uniformly on every compact set in 
$\widetilde{X}\times\widetilde{Y},$
for all multiindices $\alpha,\beta$ and all $N$ in $\mathbb{N}.$
Notice that this statement does not depend on the choice of $t,s$
and on the local coordinates $x,y.$

\par 
Our main result tells us that $K$ is negligible near every point
$(x_{0},y_{0})$ with $x_{0}\neq y_{0}$ and that for $(x,y)$ near
a diagonal point $(x_{0},x_{0})$\[
K_{t,s}(x,y)=b(x,y;\frac{1}{k})e^{k\psi(x,y)}t(x)^{k}s(y)^{*k}+R(x,y),\,\,\, 
R\,\hbox{ negligible,}\]
 where $\psi$ is smooth function with $\psi(x,x)=0,$ 
$\textrm{Re\,$\psi(x,y)$}\backsim-\left|x-y\right|^{2}$
and \[
b(x,y;\frac{1}{k})\sim k^{n}(b_{0}(x,y)+b_{1}(x,y)\frac{1}{k}+...)\]
in 
$C^{\infty}({\rm neigh}(x_{0},x_{0});
\mathcal{L}(\pi_{2}^{*}(\Lambda^{0,q}(T^{*}X\otimes 
L^{k})),\pi_{1}^{*}(\Lambda^{0,q}(T^{*}X\otimes L^{k})))$. Moreover,
let $\widetilde{\mathcal{C}}$ be the graph of $\frac{1}{i}d\psi$
in $T^{*}$$\widetilde{X}\times T^{*}\widetilde{X}$ over the diagonal.
Then $\widetilde{\mathcal{C}}$ locally represents the graph of the
canonical connection $\nabla$ of $L\otimes L^{*}$ over the diagonal
in $X\times X$ and the semiclassical wave front of $K$. See Theorem \ref{gl1} 
and the preceding explanations in Section \ref{gl} for a more precise 
local statement.

We will also explore some relations to the work \cite{s-z} of B.Shiffman
and S.Zelditch, where so called asymptotic holomorphic sections on
symplectic manifolds are studied.

\subsection{Overview}

After locally fixing a unitary frame for $L,$ we identify the Hodge
Laplacian $\Delta_{q}$ acting on $(0,q)-$forms with values in $L^{k},$
with a local semiclassical differential operator (setting $h=1/k).$
Since the curvature form of $L$ is assumed to be non-degenerate the
characteristic variety $\Sigma$ of $\Delta_{q}$ is symplectic. Modifying
the approach in \cite{MeSj1} we then construct associated local
asymptotic heat kernels in Section \ref{heat} and investigate the
limit when the time variable tends to infinity. In Section \ref{pi}
it is shown that the limit operator is an asymptotic local projection
operator. The complex canonical relation of the local projection operators
is expressed in terms of the stable outgoing and incoming manifolds
associated to $\Sigma$ in Section \ref{heat}. Assuming, in Section
\ref{gl}, that the number of negative eigenvalues of the curvature
of $L$ is equal to $q$ everywhere on $X$ we get a complete asymptotic
expansion of the global projection operator $\Pi_{q}$. In Section
\ref{ccs} we investigate some relations to \cite{s-z},
where so called asymptotic holomorphic sections on symplectic manifolds
are studied. We introduce a certain almost complex structure, closely
related to the stable manifolds introduced in Section \ref{heat},
making the curvature form of $L$ positive. It is shown that for $k$
sufficiently large the dimension of the null space of $\Delta_{q}$
coincides with the dimension of the corresponding space of asymptotically
holomorphic sections (after a suitable twisting of $L$). In Section
\ref{sec:ex} the interplay between different complex structures is
illustrated by homogeneous line bundles over flag manifolds.

\subsection{Historical remarks}
\par Most of the earlier results 
concern the positively curved case $n_-=0$. G.~Tian \cite{Ti}, followed by
W.~Ruan \cite{Ru} and Z.~Lu \cite{Lu}, computed increasingly many terms of 
the \asy{} 
expansion on the diagonal, using Tian's method of peak
solutions. T.~Bouche \cite{Bou} also got the leading term using heat
kernels.  
\par S.~Zelditch \cite{Ze}, D.~Catlin \cite{Ca} established the complete 
\asy{}
expansion 
at $x=y$ by using a result of L.~Boutet de Monvel, J.~Sj\"ostrand \cite{BoSj} 
for the
asymptotics of the Szeg\"o kernel on a strictly pseudoconvex \bdy{}
(after the pioneering work of C.~Fefferman \cite{Fe}), here on the \bdy{} 
of the unit
disc bundle, and a reduction idea of L.~Boutet de Monvel, V.~Guillemin 
\cite{BoSj}. 
Scaling asymptotics away from the diagonal was obtained later by P.~Bleher,
B.~Shiffman, S.~Zelditch \cite{BlShZe} and the full asymptotics by L.~Charles
\cite{Ch}, using again the reduction method.  In the recent work
\cite{BeBeSj} B.~Berndtsson and the authors have worked out a short
and direct proof for the full \asy{}s.

\par In more general situations, full \asy{} expansions on the 
diagonal and in some sense away from the diagonal were obtained by 
X.~Dai, K.~Liu, X.~Ma \cite{DaLiMa}. (See also \cite{MaMar} for related 
spectral results).

\par Without a positive curvature assumption there have been fewer results. 
J.M.~Bismut \cite{Bi} used the heat kernel method in his approach to 
Demailly's 
\hol{} Morse inequalities. Using local holomorphic Morse inequalities \cite{Be}, 
the leading asympotics of the Hodge projections were obtained by the first
author in \cite{b} without assuming that the curvature is non-degenerate. 
X.~Ma has pointed out to us that the method and results of \cite{DaLiMa}
can be extended to the case of non-positive \hol{} line bundles by using a
spectral gap estimate from \cite{MaMar} and this was recently carried out in
the preprint \cite{MaMar3}.  The result of Theorem \ref{gl1} was announced
in \cite{Sj3}. The precise description of the structure of the singularity 
in our result could be a first step towards the study of interesting 
tunneling phenomena in the case the curvature is of variable rank. 

\newcommand{\boldsymbol}[1]{\mbox{\boldmath $#1$}}
\medskip
\par\noindent \bf Acknowledgements. \rm  The first author has been
partially supported by a Marie Curie grant. The second author has 
benefitted from the hospitality of Chalmers and Gothenburg University in 
2000--02. We express our gratitude to Bo Berndtsson for many 
stimulating discussions and for continued interest in this work. We have 
also benefitted from discussions with L.~Boutet de Monvel, L.~Charles, 
X.~Ma, H.~Sepp\"anen and G.~Zhang, as well as to M. Shubin who suggested a 
similar problem to one of us in 1994.
\section{Holomorphic line bundles and the $\overline{\partial }$-complex, 
a review.}\label{bas}
Let $L$ be a Hermitian \hol{} line bundle over $X$. Later, we shall use 
a local \hol{} non-vanishing section $s$. We write the point-wise norm of 
$s$ as 
$$
\vert s\vert ^2=\vert s\vert _{h_1}^2=: e^{-2\phi }.$$ The curvature form of $L$ 
can be identified with the Levi form $\partial \overline{\partial }\phi $.

\par Add a Hermitian metric on $T^{1,0}X$:
\ekv{bas.6}
{
H(\nu ,\mu )=
\sum H_{j,k}\nu _k\overline{\mu }_j,
\hbox{ if }\nu =\sum \nu 
_j{\partial \over \partial z_j},\ 
\mu =\sum \mu _j{\partial \over \partial 
z_j}.}
We have a natural duality between $T^*_{1,0}X$ and $T^{1,0}X$, satisfying 
$$
\langle dz_j,{\partial \over \partial z_k}\rangle =\delta _{j,k},
$$
so if $\omega =\sum \omega _kdz_k$, then $\langle \omega ,\nu \rangle 
=\sum \omega _j\nu _j$. For each $x\in X$, we can choose $z_1,...,z_n$ centered at $x$ so 
that 
$$
H_{j,k}(x)=\delta _{j,k};\ H_x({\partial \over \partial z_j},{\partial 
\over \partial z_k})=\delta _{j,k}.
$$

The metric $H$ also determines a metric on $\Lambda ^{0,q}(T^*X)$ such that
in the special coordinates above, we have that
$$
d\overline{z}_{j_1}\wedge ...\wedge d\overline{z}_{j_q},\ 1\le 
j_1<j_2<...<j_q\le n,
$$
is an orthonormal basis of $\Lambda  ^{0,q}T_x^*X$. Then we have a 
natural 
metric also on $L\otimes \Lambda ^{0,q}T^*X$.

Let us also fix some smooth positive integration density $m(dz)$ on $X$. 
(For instance, we can take $m(dz)=\omega _n(dz)$; the induced volume 
form.) 
Then we get a natural scalar product on 
$$
{\cal E}^{0,q}(L)=C^\infty (X;L\otimes {\Lambda  }^{0,q}T^*X),
$$
so if 
$$
\overline{\partial }:\ ..\to {\cal E}^{0,q}(L)\to {\cal E}^{0,q+1}(L)\to ..
$$
is the $\overline{\partial }$ complex, then 
$$
\overline{\partial }^*:\ ..\leftarrow {\cal E}^{0,q}(L)\leftarrow {\cal 
E}^{0,q+1}(L)\leftarrow ...
$$
is also a well-defined complex.

\par If $\omega $ is a 0,1-form, let $\omega ^\rfloor 
:\Lambda ^{0,q+1}T_x^*X\to \Lambda ^{0,q}T_x^*X  $ be the adjoint of 
left 
exterior multiplication $\omega ^\wedge :\Lambda  ^{0,q}T_x^*X\to 
\Lambda ^{0,q+1}T_x^*X$. Here we use the Hermitian inner product $H^*$ 
on $\Lambda ^{0,q}T_x^*X$ that is naturally obtained from
$H$.  Without that inner product, we can still define $\nu ^\rfloor
:\Lambda ^{0,q+1}T_x^*X\to
\Lambda ^{0,q}T_x^*X$, when $\nu =\sum \nu _j{\partial \over \partial 
\overline{z}_j}$ is a \vf{} of type 0,1, as the transpose of $\nu 
^\wedge :\Lambda ^{0,q}T_xX\to\Lambda ^{0,q+1}T_xX $. We have the 
standard identity,
$$
\omega ^\wedge \nu^\rfloor+\nu ^\rfloor\omega ^\wedge =\langle \omega 
,\nu \rangle {\rm id}. 
$$
In the present case we have the analogous identity,
\ekv{bas.12}
{
\omega_1 ^\wedge \omega _2^\rfloor+\omega _2 ^\rfloor\omega_1 ^\wedge =
H^*(\omega _1,\omega _2) {\rm id}, 
}
when $\omega _1$, $\omega _2$
are (0,1)-forms. Notice also that $\omega _2^\rfloor$ depends 
anti-linearly on $\omega _2$.

\par Let $e_1(z),...,e_n(z)$ be an
\on{} frame for $\Lambda ^{0,1}T^*X$. Let $Z_1(z),...,Z_n(z)$ be the dual 
basis of $\Lambda ^{0,1}TX$, so that on scalar \fu{}s,
$$\db =\sum_1^n e_j(z)^\wedge \otimes Z_j(z,{\partial \over \partial 
\overline{z}}). 
$$

\par If $f(z)e_{j_1}\wedge ...\wedge e_{j_q}$ is a typical term in a 
general $(0,q)$-form, we get 
\begin{eqnarray*}&&\db (f(z)e_{j_1}\wedge ...\wedge e_{j_q})\\ &=&\sum 
_{j=1}^n Z_j(f)e_j^\wedge 
e_{j_1}\wedge ...\wedge e_{j_q}+\sum_{k=1}^q (-1)^{k-1} f(z)e_{j_1}\wedge 
..\wedge (\db e_{j_k})\wedge ..\wedge e_{j_q}\\
&=&(\sum_{j=1}^nZ_j(f)e_j^\wedge )e_{j_1}\wedge ...\wedge 
e_{j_k}+(\sum_{j=1}^n (\db e_j)^\wedge e_j^\rfloor )(f(z)e_{j_1}\wedge 
...\wedge e_{j_q}).
\end{eqnarray*}
So for the given \on{} frame we have the identification
\ekv{db.1}
{
\db \equiv \sum_{j=1}^n (e_j^\wedge \otimes Z_j+(\db e_j)^\wedge 
e_j^\rfloor )
}
and correspondingly
$$
\overline{\partial }^*\equiv \sum_{j=1}^n (e_j^\rfloor\otimes Z_j^*+e_j^\wedge 
(\overline{\partial }e_j)^\rfloor ),
$$
where $Z_j^*$
is the formal complex adjoint of $Z_j$ in $L^2(m)$.

\par If $s$ is a 
trivializing local \hol{} section of $L$, then $s^k$ is a 
trivializing local section of $L^k$, and the 
corresponding metric $h_k$ on $L^k$ satisfies
$$
\vert s^k\vert ^2 _{h_k}=\vert s\vert^{2k} _{h_1}=e^{-2k\phi (z)}.
$$
Hence if 
\begin{eqnarray*}
\widetilde{\omega }&=&s^k\omega \in {\cal E}^{0,q}(X;L^k),\\
\widetilde{w}&=&s^kw\in {\cal E}^{0,q+1}(X;L^k),
\end{eqnarray*}
we get for $\overline{\partial }$, $\overline{\partial }^*$, acting on 
$(0,q)$-forms with coefficients in $L^k$: 
\begin{eqnarray*}
\db (s^k\omega )&=&s^k\sum_{j=1}^n(e_j^\wedge \otimes Z_j+(\db e_j)^\wedge 
e_j^\rfloor )\omega ,\\ 
\db ^*(s^kw)&=&s^k\sum_{j=1}^n(e_j^\rfloor \otimes 
(Z_j^*+2k\overline{Z}_j(\phi ))+e_j^\wedge (\db e_j)^\rfloor )w.
\end{eqnarray*}
\par We next derive more symmetric representations for $\db$, $\db^*$ in 
spaces without exponential weights, by 
using the following local representation,
\ekv{db.2}
{
\widetilde{\omega }=(se^\phi )^k\widehat{\omega }\in {\cal E}^{0,q}(X;L^k),
}
so that 
$$
{\cal E}^{0,q}(X)\ni \widehat{\omega }\mapsto (se^\phi )^k\widehat{\omega 
}\in {\cal E}^{0,q}(X;L^k)
$$
is locally unitary in view of the fact that $\vert s(x)e^{\phi (x)}\vert 
_{h_1(x)}=1$:
\ekv{db.3}
{
\int \vert \widetilde{\omega }(x)\vert ^2_{h_k(x)\otimes H}m(dx)=\int 
\vert 
\widehat{\omega }(x)\vert _{H(x)}^2 m(dx).} 
Using \no{db.1}, which makes sense 
directly on elements of ${\cal E}^{0,q}(X,L^k)$, we get 
\ekv{db.4}
{
\db \widetilde{\omega }=(se^\phi )^k\db_s \widehat{\omega },
}
where, 
\ekv{db.5}
{
\db_s\widehat{\omega }=\sum_{j=1}^n (e_j^\wedge \otimes (Z_j+kZ_j(\phi 
))+(\db e_j)^\wedge e_j^\rfloor ).} Now the formal adjoint of $\db_s$ for 
the scalar product given by the \rhs{} of \no{db.3} is 
\ekv{db.6}
{
\db_s^*\widehat{w}=\sum_{j=1}^n (e_j^\rfloor\otimes 
(Z_j^*+k\overline{Z_j(\phi )})+e_j^\wedge (\db e_j)^\rfloor ),
} 
where in view of the unitarity of the relation \no{db.2},
\ekv{db.7}
{\db^*\widetilde{w}=(se^\phi )^k\db_s^*\widehat{w},}
where 
\ekv{db.8}
{
\widetilde{w}=(se^\phi )^k\widehat{w}.
}

\par Now rewrite things semiclassically. Put 
\ekv{db.16}
{
h={1\over k},
}
\ekv{db.17}
{
h\db_s=\sum_{j=1}^n (e_j^\wedge \otimes (hZ_j+Z_j(\phi ))+h(\db 
e_j)^\wedge e_j^\rf ),
}
\ekv{db.18}
{
h\db_s^* =\sum_{j=1}^n (e_j^\rf \otimes (hZ_j^*+\overline{Z_j(\phi 
)})+he_j^\wedge (\db e_j)^\rf ).}
Here $hZ_j$ is a semiclassical \dop{}.

\begin{prop}\label{db0}  
Using the representation \no{db.2}, we can identify the Hodge Laplacian with
\begin{eqnarray}\label{db.20}
&&\Delta=\\&&(h\db _s)(h\db _s^*)+(h\db _s^*)(h \db _s)=\nonumber\\
&&\sum_{j=1}^n 1\otimes (hZ_j^*+\overline{Z_j(\phi )})(hZ_j+Z_j(\phi 
))+\sum_{j,k}e_j^\wedge e_k^\rf \otimes [hZ_j+Z_j(\phi 
),hZ_k^*+\overline{Z_k(\phi )}]\nonumber\\
&&\hskip 3truecm +{\cal O}(h)(hZ+Z(\phi ))+{\cal 
O}(h)(hZ^*+\overline{Z(\phi )})+{\cal O}(h^2),\nonumber
\end{eqnarray}
where ${\cal O}(h)(hZ+Z(\phi ))$ indicates a remainder term of the form 
$h\sum_k a_k(z)(hZ_k+Z_k(\phi ))$ with $a_k$ smooth, matrix-valued, and 
similarly for the two other remainder terms in \no{db.20}.
\end{prop}
\par\noindent \bf Proof. \rm We make a straight forward calculation.
\begin{eqnarray*}
&&(h\db_s)(h\db_s)^*+(h\db_s)^*(h\db_s)=\\
&&\sum_{1\le j,k\le n}\Big( (e_j^\wedge \otimes (hZ_j+Z_j(\phi 
)))(e_k^\rf\otimes (hZ_k^*+\overline{Z_k(\phi )}))\\&&+(e_k^\rf\otimes 
(hZ_k^*+\overline{Z_k(\phi )}))(e_j^\wedge \otimes (hZ_j+Z_j(\phi )))\\
&&+(e_j^\wedge \otimes (hZ_j+Z_j(\phi 
)))(he_k^\wedge (\db e_k)^\rf )+(he_k^\wedge (\db e_k)^\rf )(e_j^\wedge 
\otimes (hZ_j+Z_j(\phi 
)))\\ && +h((\db e_j)^\wedge e_j^\rf )(e_k^\rf \otimes 
(hZ_k^*+\overline{Z_k(\phi )}))+(e_k^\rf \otimes 
(hZ_k^*+\overline{Z_k(\phi )}))h((\db e_j)^\wedge e_j^\rf )\\
&&+h((\db e_j)^\wedge e_j^\rf )he_k^\wedge (\db e_k)^\rf +
he_k^\wedge (\db e_k)^\rf h((\db e_j)^\wedge e_j^\rf )\Big) .
\end{eqnarray*}
Using \no{bas.12}, we see that the sum of the first two terms inside the 
general 
term of the sum is equal to \def\rf{\rfloor}
\begin{eqnarray*}
&&(e_j^\wedge e_k^\rf +e_k^\rf e_j^\wedge )\otimes 
((hZ_k^*+\overline{Z_k(\phi )})(hZ_j+Z_j(\phi )))\\
&&\hskip 4 truecm +e_j^\wedge e_k^\rf 
[hZ_j+Z_j(\phi ),hZ_k^*+\overline{Z_k(\phi )}]\\
&&= \delta _{j,k}(hZ_k^*+\overline{Z_k(\phi )})(hZ_k+Z_k(\phi 
))+e_j^\wedge 
e_k^\rf [hZ_j+Z_j(\phi ),hZ_k^*+\overline{Z_k(\phi )}].
\end{eqnarray*}
The proposition follows.\hfill{$\#$}\medskip

\par Let $q_j$ be the semiclassical principal symbol of $hZ_j+Z_j(\phi )$, 
that we shall write down more explicitly later, viewed 
as a \fu{} on the "real" cotangent space $T^*X$. (We refer to 
\cite{Ro,DiSj} for standard terminology about \scl{} \pop{}s, and to 
\cite{HeSj, SjZw} for the fact that the Weyl quantization permits to define 
the symbol of such an \op{} modulo ${\cal O}(h^2)$ even on a \mfld{}.)
The semiclassical principal 
symbol of $\Delta $ is 
\ekv{db.21}
{
p_0= 1\otimes \sum_{j=1}^n \overline{q}_jq_j.
}
The semiclassical{} subprincipal symbol of $\Delta $ is a well-defined 
endomorphism of $\Lambda  ^{0,q}T^*X$ at every point $(x,\xi )\in\Sigma $ 
on the doubly 
characteristic \mfld{} $\Sigma \subset T^*X$, given by $q_1=...=q_n=0$. 
For 
an \op{} of the form $(hZ_k^*+\overline{Z_k(\phi )})(hZ_j+Z_j(\phi ))$ 
this subprincipal symbol is given by ${h\over 2i}\{ \overline{q}_k,q_j\}$ 
and the contribution from the double sum in \no{db.20} to the subprincipal 
symbol of $\Delta $ is 
$$
{h\over i}\sum_{j,k}e_j^\wedge e_k^\rf \otimes \{ \overline{q}_j,q_k\} .
$$
Thus on $\Sigma $, we get the subprinicipal symbol of $\Delta $:
\ekv{db.22}
{
hp_1=h(1\otimes \sum_j -{1\over 2i}\{ q_j,\overline{q}_j\} 
+\sum_{j,k}e_j^\wedge e_k^\rf {1\over i}\{ q_j,\overline{q}_k\} ).
}
Since $p_1$ is invariantly defined on $\Sigma $ as well as the first sum, 
the double sum is also invariantly defined.
\par To compute further, we choose \hol{} coordinates $z_1,...,z_n$, 
$z_j=x_j+iy_j$. We make the following fiberwise bijections between 
$\Lambda  ^{1,0}T^*X$, $T^*X$, $\Lambda  ^{0,1}T^*X$:
\ekv{db.22.5}
{
\sum_1^n \zeta _jdz_j\leftrightarrow \Re (\sum_1^n \zeta 
_jdz_j)\leftrightarrow \sum_1^n \overline{\zeta }_jd\overline{z}_j.
}
Writing 
$$
\zeta _j=\xi _j-i\eta _j,
$$
we get 
$$
\Re (\sum \zeta _jdz_j)=\sum (\xi _jdx_j+\eta _jdy_j),
$$
so in local coordinates, we have bijections between 
$$
(z,\zeta )\in{\Lambda  }^{1,0}T^*X,\ (x,y;\xi ,\eta )\in T^*X,\ 
(z,\overline{\zeta })\in {\Lambda }^{0,1}T^*X .
$$
The semiclassical symbol of $h{\partial \over \partial 
\overline{z}_j}={1\over 
2}(h{\partial \over \partial x_j}+ih{\partial \over \partial y_j})$ is 
${i\over 2}(\xi _j+i\eta _j)={i\over 2}\overline{\zeta }_j$. Hence the 
symbol of 
$$
h{\partial \over \partial \overline{z}_j}+{\partial\phi  \over \partial 
\overline{z}_j}\hbox{ is }{i\over 2}\overline{\zeta }_j+{\partial \phi 
\over 
\partial \overline{z}_j},
$$
so in the coordinates $(z,\zeta )$, the \e{} for $\Sigma $ becomes:
$$
\overline{\zeta }_j=-{2\over i}{\partial \phi \over \partial 
\overline{z}_j},
$$
or equivalently,
\ekv{db.23}
{
\zeta _j={2\over i}{\partial \phi \over \partial z_j},\ j=1,2,..,n.
}

\par For later use we here compute the principal symbol $q_j$ of 
$hZ_j+Z_j(\phi )$: Let the \on{} frame $e_1,...,e_n$ be given by 
$$
e_j(z)=\sum_k a_{j,k}(z)d\overline{z}_k,
$$  
and the corresponding dual basis $Z_1,...,Z_n$ of $\Lambda ^{0,1}T_z^*X$ 
by
$$
Z_j=\sum_k b_{j,k}{\partial \over \partial \overline{z}_k},
$$
where the invertible matrices $(a_{j,k})$ and $(b_{j,k})$ are related by
$$
\trans (b_{j,k})(a_{j,k})=1.
$$
Then it follows from the calculations above that 
\ekv{db.23.5}
{q_j=\sum_k b_{j,k}({i\over 2}\overline{\zeta }_k+{\partial \phi \over 
\partial \overline{z}_k}).}
\begin{prop}\label{db0.5} In the $(z,\zeta )$-coordinates, the Poisson 
bracket $\{ f,g\}$ of two $C^1$-\fu{}s $f,g$ is given by
\ekv{db.25}
{
{1\over 2}\{ f,g\} ={1\over 2}H_fg=({\partial f\over \partial \zeta }\cdot 
{\partial g\over  \partial z}+{\partial f\over \partial 
\overline{\zeta } }\cdot {\partial g\over  \partial \overline{z}})-  
( {\partial 
f\over \partial z }\cdot {\partial g\over  \partial \zeta } 
+{\partial f\over \partial \overline{z} }\cdot {\partial g\over  \partial 
\overline{\zeta }})
}
\end{prop}
\par\noindent \bf Proof. \rm
Consider the real canonical 
1-form on $T^*X$:
$$
\Re (\sum \zeta _jdz_j)=\sum (\xi _jdx_j+\eta _jdy_j).
$$
Hence the real symplectic form becomes 
$$
d(\sum (\xi _jdx_j+\eta _jdy_j))=\Re (\sum d\zeta _j\wedge dz_j)=\Re 
\sigma =:\omega ,
$$
where $\sigma =\sum d\zeta _j\wedge dz_j$. If $f$
is a smooth real \fu{} on the real phase space, the corresponding Hamilton 
field $H_f$ is given by 
\ekv{db.24}
{
\langle \omega ,t\wedge H_f\rangle =\langle t,df\rangle .
}
With $t=2\Re \sum (a_j{\partial \over \partial z_j}+b_j{\partial \over 
d\zeta _j})$, the \rhs{} becomes 
$$2\Re \sum (a_j{\partial f\over \partial z_j}+b_j{\partial f \over 
\partial 
\zeta _j}),$$
while the \lhs{} is equal to 
$$
\Re \langle \sigma ,t\wedge H_f\rangle =\Re \sum (b_j\langle 
dz_j,H_f\rangle -a_j\langle d\zeta _j,H_f\rangle ).
$$
Varying $t$, we conclude that 
$$\langle dz_j,H_f\rangle =2{\partial f\over 
\partial \zeta _j},\ \langle d\zeta _j,H_f\rangle =-2{\partial f\over 
dz_j},
$$
so 
$$
{1\over 2}H_f=({\partial f\over \partial \zeta }\cdot {\partial\over  
\partial
z}-{\partial  f\over \partial z }\cdot {\partial\over  \partial \zeta }) + 
( {\partial f\over \partial 
\overline{\zeta } }\cdot {\partial\over  \partial \overline{z}} 
-{\partial f\over \partial \overline{z} }\cdot {\partial\over  \partial 
\overline{\zeta }}).
$$
In particular, we get \no{db.25}
This expression now extends to the case when $f,g$ are complex-valued 
\fu{}s which completes the proof. \hfill{$\#$}\medskip

\par Of course \no{db.25} can also be obtained by straight forward 
calculation from 
\ekv{db.26}
{
\{ f,g\} ={\partial f\over \partial\xi  }
{\partial g\over \partial x }+
{\partial f \over \partial \eta  }{\partial g \over \partial y}
-{\partial f\over \partial x }{\partial g\over \partial \xi }
-{\partial f \over \partial y }{\partial g \over \partial \eta },\ 
{\partial \over \partial x}={\partial \over \partial z }+
{\partial \over \partial \overline{z} },\,
{\partial \over \partial y}={1\over i}({\partial \over 
\partial\overline{z} }
-{\partial \over \partial z}), ...
}

\par Now return to the expressions \no{db.20}, \no{db.21}.  If $z_0$
is a fixed point, we choose \hol{} coordinates $z_1,...,z_n$ as above in
such a way that $Z_j={\partial \over \partial \overline{z}_j}$, 
$e_j=d\overline{z}_j$
at $z_0$.  Then $b_{j,k}(z_0)=\delta _{j,k}$ in \no{db.23.5} and at the
corresponding point $\rho _0=(z_0,\zeta _0)\in
\Sigma $, we have 
$$
\{ q_j,\overline{q}_k\} (\rho _0)=\{ {i\over 2}\overline{\zeta 
}_j+{\partial \phi \over \partial \overline{z}_j}, -{i\over 2}\zeta 
_k+{\partial \phi \over \partial z_k}\}.
$$
Applying \no{db.25}, we now get
$$
{1\over 2}\{ q_j,\overline{q}_k\}={i\over 2}{\partial ^2\phi \over 
\partial 
\overline{z}_j\partial z_k}+{\partial ^2\phi\over  \partial z_k\partial 
\overline{z}_j }{i\over 2}=i{\partial ^2\phi \over \partial 
\overline{z}_j\partial z_k}.
$$
We rewrite this as 
\ekv{db.27}
{
{1\over 2i}\{ {q}_j,\overline{q}_k\} ={\partial ^2\phi \over \partial 
\overline{z}_j\partial z_k},
}
and recognize here the coefficients of the Levi-matrix appearing also in 
$\overline{\partial }\partial \phi $.
\begin{Prop}\label{db2} $\Sigma $ is symplectic at a point 
$(z_0;\xi _0,\eta _0)$ iff $({\partial ^2\phi \over \partial 
\overline{z}_j\partial z_k})(z_0)$ is \nondeg{}. Indeed, if we identify 
$\Lambda  ^{1,0}T^*X$ and $T^*X$, by means of the first bijection in 
\no{db.22.5}, 
then the real symplectic form $\omega $ becomes $\Re (\sum d\zeta _j\wedge 
dz_j)$
and its restriction to $\Sigma $ can be identified with ${2\over 
i}\overline{\partial }\partial \phi $.\end{Prop}
\par\noindent \bf Proof. \rm With the above mentioned identification, 
$\Sigma $ takes the form \no{db.23} which can be written more invariantly 
as 
\ekv{db.34}
{
\zeta \cdot dz={2\over i}\partial \phi .
}
Hence,
$$
{\sigma _\vert}_{\Sigma }=d\sum_1^n {2\over i}{\partial \phi \over 
\partial z_j}\wedge dz_j=\sum_{j=1}^n\sum_{k=1}^n {2\over i}{\partial 
^2\phi \over \partial \overline{z}_k\partial z_j}d\overline{z}_k\wedge 
dz_j={2\over i}\overline{\partial }\partial \phi .
$$ 
This is a real form, so it is also the restriction to $\Sigma $ of $\Re 
\sigma $ and it is \nondeg{} precisely when $({\partial ^2\phi \over 
\partial \overline{z}_j\partial z_k})$ is (cf 
\cite{Sj1}).\hfill{$\#$}\medskip

\par Back to the general case, we recall the condition for having the 
apriori estimate
\ekv{db.28'}
{h\Vert u\Vert +\sum \Vert (hZ_j+Z_j(\phi ))u\Vert + \sum \Vert 
(hZ_j^*+\overline{Z_j(\phi )})u\Vert \le C\Vert \Delta _qu\Vert ,}
for $u\in 
C_0^\infty ({\rm neigh\,}(z_0);{\Lambda  }^{0,q}T^*X).$
\begin{Prop}\label{db1} \no{db.28'} does not hold precisely when $n_-\le 
q\le 
n-n_+$, where $(n_+,n_-)$ is the signature of $({\partial ^2\phi \over 
\partial 
\overline{z}_j\partial z_k}(z_0))$.\end{Prop}

This is essentially well-known since the $\db$-estimates of L.~H\"ormander
(see \cite{Ho}), and in the context of more general hypoelliptic operators it
was obtained in \cite{Sj2} in the \nondeg{} symplectic case.  The
result will not be used explicitly since the heat \e{} method below will
give enough control (and would allow to recover it easily).

\section{The associated heat \e{}s.}\label{heat}
\par We work locally near a point $z_0\in X$, where 
\ekv{heat.1}
{
({\partial ^2\phi \over \partial \overline{z}_j\partial z_k})\hbox{ is 
\nondeg{} of signature }(n_+,n_-),
}
so that the characteristic \mfld{} $\Sigma $ of $\Delta _q$ is symplectic. 
We review some results of A.~Menikoff, J.~Sj\"ostrand \cite{MeSj1}, \cite{MeSj2} 
that apply to the 
present situation with minor changes:

\par In those works, we considered  a scalar classical \pop{} with 
principal symbol $p_0$ vanishing to precisely the second order on a conic
symplectic sub\mfld{} of $T^*X$.  In the present work, we have a
semiclassical{} \dop{} with a leading symbol $p_0$ in \no{db.1} that we can
view as scalar; $p_0=\sum_1^n\overline{q}_jq_j$ and $p_0$ is no longer
homogeneous, and $\Sigma $ is no longer conic in the fiber variables.

\par In this section we consider the \pb{}:
\ekv{heat.2}
{
(h\partial _t+\Delta _q)u(t,x)=0,\ u(0,x)=v(x).
}
We shall apply the standard WKB construction of an approximative solution 
\op{} and apply arguments from \cite{MeSj1} together with a "Witten trick" to 
get additional properties to be used later. See Proposition \ref{heat1.1} 
for the precise statement about the solution to \no{heat.2}.

\par We forget about  most of the complex structure of $X$ and work in some 
smooth local coordinates $x=(x_1,...,x_{2n})$ defined on 
$\widetilde{X}\subset\subset X$. At least for small $t\ge 
0$, we look for an approximate solution of \no{heat.2} of the form 
$u(t,x)=U(t)v(x)$, 
\ekv{heat.3}
{
u(t,x)={1\over (2\pi h)^{2n}}\iint e^{{i\over h}(\psi (t,x,\eta )-y\cdot 
\eta )}a(t,x,\eta ;h)u(y)dyd\eta ,
}
where $a$ is a matrix-valued classical symbol of order $0$:
\ekv{heat.4}
{
a(t,x,\eta ;h)\sim \sum_0^\infty  a_k(t,x,\eta )h^k,\ {a_\vert}_{t=0}=1,
}
and $\psi $ with $\Im \psi \ge 0$ should solve the eikonal \e{},
\ekv{heat.5}
{
i\partial _t\psi (t,x,\eta )+p_0(x,\psi '_x(t,x,\eta ))=0+{\cal O}((\Im 
\psi )^\infty ),\ {\psi _\vert}_{t=0}=x\cdot \eta .
} The amplitude $a$ is determined by a sequence of transport \e{}s that
will be reviewed later. (Here we follow the convention that $u={\cal O}((\Im 
\psi )^\infty )$ means that $u={\cal O}((\Im 
\psi )^N)$ for every $N\ge 0$, \ufly{} or locally \ufly{} depending on the 
context.)

\par According to the general theory in \cite{MelSj,MelSj2}, this equation 
can be solved 
locally, provided that we also denote by $p_0$ an almost \hol{} extension. 
The general theory also tells us that $U(t)$ is associated to a \ctf{},
\ekv{heat.6}
{
\kappa _t=\exp (-itH_{p_0}).
}
(Here $\kappa _t$ depends slightly on the choice of almost \hol{} 
extension 
of $p_0$, so $\kappa _t(\rho )$ is well-defined only up to $\vert \Im \rho 
\vert ^\infty $. In \cite{MelSj2} we also made the assumption that 
$p_0(x,\xi )$ is positively homogeneous of degree 1 in $\xi $, but as 
noticed for instance in \cite{MeSj1} and will be reviewed in the proof of 
Proposition \ref{heat1.1}, one can easily reduce the general 
case to the homogeneous one, by adding a variable $x_0$ and consider the 
homogeneous symbol $\xi _0 p(x,\xi /\xi _0)$, then restrict the results 
to $\xi _0=1$.)

\par So far, we only used the non-negativity of (the real part of) $p_0$. Now 
we use that $p_0\sim {\rm dist\,}(\cdot ,\Sigma )^2$. It follows that 
\ekv{heat.7}
{
\psi (t,x,\eta )=x\cdot \eta +{\cal O}(t\,{\rm dist\,}(x,\eta ;\Sigma )^2),
}
\ekv{heat.8}
{
\Im \psi (t,x,\eta )\sim t\,{\rm dist\,}(x,\eta ;\Sigma )^2,
}
for $0\le t\le t_0$, and $t_0>0$ fixed. Correspondingly, we have 
\ekv{heat.9}
{
{{\kappa _t}_\vert}_{\Sigma }={\rm id},
}
\begin{eqnarray}
\label{heat.10}
&&\hbox{When }t>0,\ \kappa _t\hbox{ is a \st{} positive canonical}
\\ &&\hbox{transformation with } {\rm graph\,}(\kappa _t)\cap 
(T^*X)^2={\rm diag\,}(\Sigma
\times \Sigma ).\nonumber\end{eqnarray} 
Recall that a positive \ctf{} is strictly positive if the
${\rm graph\,}\kappa $ intersects $T^*X\times T^*X$ cleanly along a smooth
sub\mfld{}.  Thanks to these simplifying features, all essential properties
of $\psi $ and $\kappa _t$ are captured by their Taylor expansions at $t=0$
and at $\Sigma $.

\par In \cite{MeSj1} it was shown that \no{heat.5} can be solved for all 
$t\ge 0$, 
and that we have,
\ekv{heat.11}
{
\Im \psi (t,x,\eta )\sim {\rm dist\,}(x,\eta ;\Sigma )^2,
}
\ufly{} for $t\ge 1$, that \no{heat.9}, \no{heat.10} remain valid for all 
$t>0$, and 
finally that there exists a smooth \fu{} $\psi (\infty ,x,\eta )$, 
well-defined mod ${\cal O}({\rm dist\,}(x,\eta ;\Sigma )^\infty )$ such 
that for all $k,\alpha $:
\ekv{heat.12}
{
\partial _t^k\partial _{x,\eta }^\alpha (\psi (t,x,\eta )-\psi (\infty 
,x,\eta ))={\cal O}(e^{-t/C}),
}
\ufly{} on $[0,+\infty [\times \Sigma $. (In \cite{MeSj2} we also established \asy{} 
expansions 
when $t\to \infty $ in terms of exponentials in $t$.  We do not need those
improved results here.)  Here we have locally \ufly{} on 
$\widetilde{X}\times {\bf R}^{2n}$:
\ekv{heat.13}
{
\psi ( \infty ,x,\eta )=x\cdot \eta +{\cal O}({\rm dist\,}(x,\eta ;\Sigma 
)^2),\ \Im\psi ( \infty ,x,\eta )\sim {\rm dist\,}(x,\eta ;\Sigma 
)^2.
} 
Further, \def\cre{canonical relation} the \cre{} $C_\infty $ generated by 
the phase $\psi (\infty ,x,\eta )-y\cdot \eta $ is \st{} positive with 
\ekv{heat.14}
{
C_\infty \cap (T^*X\times T^*X)={\rm diag\,}(\Sigma \times \Sigma ),
}
and $C_\infty $
can be described in the following way:

\par There are two almost 
\hol{} \mfld{}s $J_+,J_-\subset T^*X^{\bf C}$ (where the latter set is the 
almost complexification of $T^*X$) intersecting $T^*X$ cleanly along 
$\Sigma $, with the following properties:
\begin{eqnarray}\label{heat.14.5} &&{\rm codim}_{\bf C}J_\pm =n,\ 
J_{\pm}\subset p_0^{-1}(0),\\
&&J_{\pm} \hbox{ are involutive and } J_-=\overline{J}_+,\nonumber\\
&&{1\over i}\sigma (t,\overline{t})>0,\, \forall t\in T_\rho 
(J_+)\setminus T_\rho (\Sigma ^{\bf C}),\,\, \rho \in \Sigma .\nonumber
\end{eqnarray}
\par Here the involutivity of $J_+$ (and similarly for $J_-$) means that 
$J_+$ is given by the equations $\widetilde{q}_1=...=\widetilde{q}_n=0$, 
where $d\widetilde{q}_1,...,d\widetilde{q}_n$ are ${\bf C}$-linearly 
\indep{} and $\{ \widetilde{q}_j,\widetilde{q}_k\} =0$ on $J_+$. Further 
the complexification $\Sigma ^{\bf C}$ is contained in $J_{+}$ and 
$H_{\widetilde{q}_1},...,H_{\widetilde{q}_n}$ span $T_\rho J_+/T_\rho 
\Sigma ^{\bf C}$. The positivity property above is equivalent to 
the fact that the Hermitian matrix $({1\over i}\{ 
\widetilde{q}_j,\overline{\widetilde{q}}_k\} )$ is positive definite. In 
terms of $J_{\pm}$, we can describe the limiting canonical relation 
$C_\infty $ 
as
$\{ (\rho ,\mu )\in J_+\times J_-;\,$ the
$n$-dimensional bicharacteristic leaves through $\rho $, $\mu $ of $J_+$ 
and $J_-$ respectively, 
intersect $\Sigma _+^{\bf C}$ at the same point $\}$.

Finally we can also view $C_\infty $ as the limit of $C_t={\rm 
graph\,}(\kappa _t)$, when $t\to +\infty $, where the convergence is 
exponentially fast (in the sense of Taylor expansions at ${\rm 
diag\,}(\Sigma \times \Sigma )$). We can also view $J_+$, $J_-$
as the stable outgoing and incoming \mfld{}s respectively, for the 
$H_{-ip}$-flow, near the fixed point set $\Sigma ^{\bf C}$. Let us also 
add 
that $J_\pm$
are uniquely determined and that in the case $n_+=n$, we can take 
$\widetilde{q}_j=q_j$.

\par Next we consider the behaviour of $a$ in \no{heat.3}, \no{heat.4}, 
where we recall 
that $a_0$, $a_1$, ... are successively determined by a sequence of 
transport \e{}s. Following \cite{MeSj1} this can be done in the following 
way, where we take some advantage of the fact that we work in the Weyl 
quantization. (See also appendix b of \cite{HeSj}.): Formally, with 
$\psi =\psi (t,\cdot ,\eta )$, $P=\Delta _q$ and with the exponent $w$ 
indicating that we take the $h$-Weyl quantization, we get 
\begin{eqnarray*}
&&e^{-i\psi}\circ P\circ e^{i\psi /h}=P(x,\psi_x '(x)+\xi ;h)^w+{\cal 
O}(h^2)=\\&&
 p(x,\psi '_x)+hp_1(x,\psi '_x)+{1\over 2}(hD_x\circ p'_\xi (x,\psi 
'_x)+p'_\xi (x,\psi '_x)\circ hD_x)+{\cal O}(h^2)= \\&&
p(x,\psi '_x)+hp_1(x,\psi '_x)+{h\over i}p'_\xi (x,\psi '_x)\cdot 
{\partial \over \partial x}+{h\over 2i}{\rm div\,}(p'_\xi (x,\psi 
'_x)\cdot {\partial \over \partial x})+{\cal O}(h^2),
\end{eqnarray*}
where the "${\cal O}(h^2)$"
refers to the action on symbols and $p_1$ is the subprincipal symbol. This 
gives the first transport \e{} for $a_0$:
$$
(\nu +{1\over 2}{\rm div\,}(\nu )+p_1)a_0=0,
$$
where
$$
\nu ={\partial \over \partial t}-ip'_\xi (x,\psi '_x)\cdot {\partial 
\over \partial x}.
$$
The higher transport equations for $a_j$, $j\ge 1$, are of the form 
$$
\nu (a_j)=F_j(t,x,a_0,...,a_{j-1}).
$$ 
Then if $a(t,x,\eta ;h)\sim\sum_0^\infty a_j(t,x,\eta )h^j$ in $C^\infty 
([0,+\infty [\times \widetilde{X}\times {\bf R}^{2n})$, we have 
$$
(h\partial _t+\Delta _q)(e^{{i\over h}\psi (t,x,\eta )}a(t,x,\eta 
;h))={\cal O}(h^\infty )
$$
locally \ufly{} on $[0,+\infty [\times \widetilde{X}\times {\bf R}^{2n}$ and 
similarly for the derivatives.
\par The discussion on page 69 in \cite{MeSj1} shows that ${\rm div\,}(\nu
)\to {1\over 2}\widetilde{{\rm tr}}\,F$ exponentially fast on $\Sigma $,
where $\widetilde{{\rm tr}}\, F=\sum f_j$, and $F$ is the fundamental
matrix of $p$ ie the linearization of $H_p$ at the point of $\Sigma $ and
has the spectrum $\sigma (F)=\{ \pm if_j\}$, $f_j\ge 0$.  In the further
discussion of the transport equations the only new feature is that $p_1$ is
now a square matrix rather than a scalar, and whenever we needed a lower
bound on $\Re p_1$, we now need a lower bound on the set of real parts of
the \ev{}s of $p_1$.  Proposition 2.2 in \cite{MeSj1} becomes
\begin{Prop}\label{heat'1}
Let $\lambda \in C(\Sigma ;{\bf R})$ satisfy 
$$
\lambda (x,\eta )<{1\over 2}\widetilde{{\rm tr}}\,F(x,\eta )+{\rm 
inf\,}\Re \sigma (p_1(x,\eta )),\ (x,\eta )\in \Sigma .
$$
Then for every compact set $K\subset \Sigma $, $j\in{\bf N}$
and $(\gamma ,\alpha ,\beta )\in{\bf N}^{1+2n+2n}$, we have 
$$
\vert \partial _t^\gamma \partial _x^\alpha \partial _\eta ^\beta 
a_j(t,x,\eta )\vert \le C_{j,\alpha ,\beta ,\gamma }e^{-t\lambda 
(x,\eta )},\ (x,\eta )\in K,\ t\ge 0.
$$
\end{Prop}

\par We are therefore interested in whether 
\ekv{heat'.1}
{
{1\over 2}\widetilde{{\rm tr}}\, F+\inf \Re \sigma (p_1)>0\hbox{ on 
}\Sigma 
}
or not. Now 
$$
p=\sum_1^n \overline{q}_jq_j,\quad H_p=\sum 
(\overline{q}_jH_{q_j}+q_jH_{\overline{q}_j}).
$$
At a given point $\rho _0\in \Sigma $, we choose the basis 
$H_{q_1},...,H_{q_n},H_{\overline{q}_1},...,H_{\overline{q}_n}$ for 
$T_{\rho _0}(T^*X)^{\bf C}/\Sigma ^{\bf C}$, and compute the linearization 
of $H_p$:
$$
H_p(\rho _0+\sum t_kH_{q_k}+\sum s_kH_{\overline{q}_k})={\cal 
O}((t,s)^2)+\sum_{j,k}t_k\{ q_k,\overline{q}_j\} H_{q_j}+\sum_{j,k}s_k\{ 
\overline{q}_k,q_j\} H_{\overline{q}_j}.
$$
So the matrix $F_p$ of the linearization is expressed in the basis above 
by
$$
{1\over i}F_p=\pmatrix{({1\over i}\{ q_k,\overline{q}_j\}) & 0\cr 0 &
({1\over i}\{ \overline{q}_k,q_j\})},
$$
where we recall \no{db.27}. Let $\mu _1,...,\mu _n$ be the \ev{}s of 
$(\partial _{\overline{z}_j}\partial _{z_k}\phi )$, with $\mu _j>0$ for 
$1\le j\le n_+$ and $\mu _j<0$ for $n_++1\le j\le n$. Then $$(i^{-1}\{ 
q_k,\overline{q}_j\} )={}\trans (i^{-1}\{ q_j,\overline{q}_k\} )\hbox{ has 
the 
\ev{}s } 2\mu _1,...,2\mu _n,$$ and $$(i^{-1}\{
\overline{q}_k,{q}_j\} )=-{}\trans (i^{-1}\{ q_j,\overline{q}_k\} )\hbox{ 
has the 
\ev{}s }-2\mu _1,...,-2\mu _n.$$
Hence the non-vanishing \ev{}s of $F_p$ are $\pm 2i\mu _1,...,\pm 2i\mu 
_n$, and 
\ekv{heat'.2}
{
{1\over 2}\widetilde{{\rm tr}}\, F_p=\mu _1+...+\mu _{n_+}-\mu 
_{n_++1}-...-\mu _n.
}

\par For the first term in \no{db.22}, we get
\ekv{heat'.3}
{
\sum_j -{1\over 2i}\{ q_j,\overline{q}_j\}=-{1\over 2i}{\rm tr\,}(\{ 
q_j,\overline{q}_k\} )=-\sum_1^n \mu _j.
}
We can also compute the \ev{}s of the matrix part of the subprincipal 
symbol 
appearing in \no{db.22} and in the subsequent remark about invariance. We 
choose \hol{} coordinates such that at the given point $z_0$: 
$Z_j=\partial  
_{\overline{z}_j}$, $e_j=d\overline{z}_j$ and moreover $(i^{-1}\{ 
q_j,\overline{q}_k\} )$ is diagonalized, equal to 
$$
\pmatrix{2\mu _1 &0 &.. &0\cr
0 &2\mu _2 &..&0\cr
.. & &..  &..\cr
0&0&..&2\mu _n}.
$$
Then 
$$
\sum_{j,k}{1\over i}\{ q_j,\overline{q}_k\} e_j^\wedge e_k^\rf =\sum_j 2\mu 
_j e_j^\wedge e_j^\rf . 
$$
On $(0,q)$-forms, the \ev{}s are the numbers 
$$
2(\mu _{j_1}+\mu _{j_2}+...+\mu _{j_q}),\hbox{ for }1\le 
j_1<j_2<...<j_q\le n.
$$
From \no{heat'.2}, \no{db.22} and the other calculations we get
$$
p_1+{1\over 2}\widetilde{{\rm tr}}\, F=-2\sum_{n_++1}^n \mu 
_j+\sum_{j,k}{1\over i}\{ q_j,\overline{q}_k\} e_j^\wedge e_k^\rf ,
$$
which on the space of $(0,q)$-forms has the \ev{}s
$$
-2\sum_{n_++1}^n \mu _j +2(\mu _{j_1}+..+\mu _{j_q}),\ 1\le j_1<...<j_q\le 
n.
$$
We see that on $\Sigma $
\ekv{heat'.4}
{
\inf \sigma (p_1+{1\over 2}\widetilde{{\rm tr}}\,F)\cases{=0,\ q=n_-\cr 
> 0,\ q\ne n_-}.
} This is the answer to the question \no{heat'.1} and Proposition
\ref{heat'1} then shows that when $q\ne n_-$, there exists a constant
$C>0$ such that
\ekv{heat.15}
{
\vert \partial _t^k\partial _{x,\eta }^\alpha a_j(t,x,\eta )\vert \le 
C_{k,\alpha 
,j}e^{-t/C},\ t\ge 0,\,(x,\eta )\in \Sigma ,
}
while in the case $q=n_-$, we have for every $\epsilon >0$:
\ekv{heat.16}
{
\vert \partial _t^k\partial _{x,\eta }^\alpha a_j(t,x,\eta )\vert \le 
C_{k,\alpha ,j,\epsilon }e^{\epsilon t},\ t\ge 0,\, (x,\eta )\in \Sigma .
}
We also notice from \cite{MeSj1}, that \no{heat.15} and 
\no{heat.16} respectively hold 
also when the initial condition in \no{heat.4} is replaced by 
${a_\vert}_{t=0}=b$ for any classical symbol $b(x,\eta 
;h)\sim\sum_0^\infty 
b_j(x,\eta )h^j$.

\par Using the particular structure of the \pb{}, we will next show
\begin{Prop} \label{heat1} Consider the case $q=n_-$
and let $a$ be the symbol in \no{heat.3}, \no{heat.4}. Then there exist 
$C>0$ and a 
classical symbol 
$$
a^\infty (x,\eta ;h)\sim \sum_0^\infty  a_j^\infty (x,\eta )h^j,
$$ 
such that 
\ekv{heat.17}
{
\vert \partial _t^k\partial _{x,\eta }^\alpha (a_j(t,x,\eta )-a_j^\infty 
(x,\eta ))\vert \le C_{k,\alpha ,j}e^{-t/C},\ t\ge 0,\, (x,\eta 
)\in\Sigma .}\end{Prop}

\par\noindent \bf Proof. \rm $a$ is determined by the initial condition 
in \no{heat.4} and the fact that 
\ekv{heat.18}
{
(h\partial _t+\Delta _q)(e^{{i\over h}\psi (t,x,\eta )}a(t,x,\eta 
;h))={\cal O}(h^\infty ),
}
locally \ufly{} in $t$, and similarly for the derivatives. Let $Z_\phi
:=h\db_s$ be given in \no{db.17} , so  that $Z^*_\phi $ is given by 
\no{db.18}. Then we
have the intertwining  properties,
\ekv{heat.19}
{
\Delta _{q+1}Z_\phi =Z_\phi \Delta _q,\ \Delta _{q-1}Z_\phi ^*=Z_\phi 
^*\Delta _q.
}
Combining this with \no{heat.18}, we get 
\ekv{heat.20}
{
(h\partial _t+\Delta _{q-1})(Z_\phi ^*(e^{{i\over h}\psi }a))={\cal 
O}(h^\infty ),
}
\ekv{heat.21}
{
(h\partial _t+\Delta _{q+1})(Z_\phi (e^{{i\over h}\psi }a))={\cal 
O}(h^\infty ).
}
Now
\ekv{heat.21.5}
{
Z_\phi^* (e^{{i\over h}\psi }a)=e^{{i\over h}\psi }\widetilde{a},\ Z_\phi 
(e^{{i\over h}\psi }a)=e^{{i\over h}\psi }\widehat{a},
}
where $\widetilde{a}$, $\widehat{a}$
are classical symbols of order 0 in $h$, and combining this with 
\no{heat.20}, 
\no{heat.21}, we see that \no{heat.15} applies to $\widetilde{a}$, 
$\widehat{a}$. 
Now, $\Delta _q=Z_\phi^* Z_\phi +Z_\phi Z_\phi^* $, so 
\ekv{heat.22}
{
\Delta _q(e^{{i\over h}\psi }a)=e^{{i\over h}\psi }b
}
where $b\sim\sum_{0}^\infty b_j(t,x,\eta )h^j$ and the $b_j$ satisfy 
\no{heat.15}.

\par Combining this with \no{heat.18}, we see that 
\ekv{heat.23}
{
h\partial _t(e^{{i\over h}\psi }a)=e^{{i\over h}\psi }c,
}
where $c$ ($=-b+{\cal O}(h^\infty )$) has the same properties as $b$. But 
$$
c=h\partial 
_ta+i(\partial _t\psi )a,
$$
so if we combine \no{heat.12}, \no{heat.16} with the fact that $c$ 
satisfies \no{heat.15}, we get 
\ekv{heat.24}
{
\vert \partial _t^k\partial _{x,\eta }^\alpha \partial _ta_j(t,x,\eta 
)\vert \le C_{k,\alpha ,j}e^{-t/C},\ t\ge 0,\, (x,\eta )\in\Sigma .
}
From this we get \no{heat.17}.\hfill{$\#$}\medskip

\par We introduce the \scl{} Sobolev space
$$
H^s({\bf R}^{2n})=\{ u\in{\cal S}'({\bf R}^{2n});\, \langle hD_x\rangle 
^su\in L^2\},\ s\in{\bf R},
$$
with the $h$-dependent norm $\Vert u\Vert _{H^s}=\Vert \langle hD_x\rangle 
^su\Vert $. Here $\langle hD_x\rangle =(1+(hD)^2)^{1/2}$. From this, we 
form $H^s_{\rm comp}(X)$, $H^s_{\rm loc}(X)$ in the usual way, when $X$
is a smooth paracompact \mfld{}, as well as $H^s(X)$, when $X$
is compact. On the Fr\'echet space $H^s_{\rm loc}(X)$, we have natural 
$h$-dependent semi-norms, so it makes sense to say that $u=u_h$ is ${\cal 
O}(h^{N_0})$ in $H^s_{\rm loc}$.

\par We now return to our local coordinate patch $\widetilde{X}\subset X$, 
and define 
\ekv{heat.24.1}
{
U(t)u(x)={1\over (2\pi h)^{2n}}\iint e^{{i\over h}(\psi (t,x,\eta )-y\cdot 
\eta 
)}a(t,x,\eta ;h) u(y)dyd\eta ,
}
with $\psi$, $a\sim\sum_0^\infty a_j(t,x,\eta )h^j$ constructed as above. 
More precisely, we can choose $a,a_j\in C^\infty ([0,\infty [\times 
\widetilde{X}\times {\bf R}^{2n})$ with the following properties:
\ekv{heat.24.2}
{
\partial _t^k\partial _x^\alpha \partial _\eta ^\beta  a_j=\cases{
{\cal 
O}_{j,\alpha ,\beta ,K}(1)e^{-t/C}, q\ne n_-\cr
{\cal 
O}_{j,\alpha ,\beta ,K,\epsilon }(1)e^{\epsilon t},\ q=n_- 
},\ (x,\eta )\in K\subset\subset \widetilde{X}\times {\bf R}^{2n},\ 
\epsilon >0,
}
\begin{eqnarray}\label{heat.24.3}
&\partial _t^k\partial _x^\alpha \partial _\eta ^\beta (a-\sum_0^{N-1}  
h^ja_j)=h^N\cases{
{\cal 
O}_{k,\alpha ,\beta ,K}(1)e^{-t/C}, q\ne n_-\cr
{\cal 
O}_{k,\alpha ,\beta ,K,\epsilon }(1)e^{\epsilon t},\ q=n_- 
},&\\ &(x,\eta )\in K\subset\subset \widetilde{X}\times {\bf R}^{2n},\ 
\epsilon >0.&\nonumber
\end{eqnarray}
Moreover, in the case when $q=n_-$, we have $a(\infty ,x,\eta 
;h)\sim\sum_0^\infty  a_j(\infty ,x,\eta )h^j$ in $C^\infty 
(\widetilde{X}\times {\bf R}^{2n})$, such that
\ekv{heat.24.4}
{
\partial _t^k\partial _x^\alpha \partial _\eta ^\beta (a_j(t,x,\eta 
)-a_j(\infty ,x,\eta ))={\cal O}_{k,\alpha ,\beta ,K}(1)e^{-t/C}
}
and similarly for $a(t,x,\eta ;h)-a(\infty ,x,\eta ;h)$. We also arrange 
so that 
\ekv{heat.24.5}{
a(0,x,\eta ;h)=1.
}

\par The construction of $\psi ,a$ can be extended in the natural way to 
the elliptic region $\vert \eta \vert \gg 1$, and here it all boils down to
Taylor expanding in $t$.  We quickly review a way of treating
this standard heat evolution problem by a simple dilation argument.  (The
reader may skip this and go directly to Proposition \ref{heat1.1}.) If $\Delta 
_q=P(x,hD_x;h)$ (say with $P(x,\xi ;h)$
denoting the Weyl symbol for our local coordinates) then in the problem 
\no{heat.2}, we let $\lambda \gg 1$ and make the change of time variable 
$s=\lambda t$, so that $\lambda ^{-1}\partial _t=\partial _s$. Then 
dividing \no{heat.2} by $\lambda ^2$, we get the new evolution \e{} 
\ekv{heat.24.5.1}
{
(\widetilde{h}\partial _s +\widetilde{P}(x,\widetilde{h}D_x,{1\over 
\lambda };\widetilde{h}))u=0,\ \widetilde{h}=h/\lambda ,
}
where
\ekv{heat.24.5.2}
{
\widetilde{P}(x,\xi ,{1\over \lambda };\widetilde{h})={1\over \lambda 
^2}P(x,\lambda \xi ;h).
}
Recall here that $P(x,\xi ;h)=p(x,\xi )+hp_1(x,\xi )+h^2p_2(x)$, where 
$p$, $p_1$, $p_2$ are polynomials in $\xi $
of degree $2$, $1$ and $0$
respectively. If we decompose into homogeneous polynomials, 
\begin{eqnarray*}
p(x,\xi )&=&p^2(x,\xi )+p^1(x,\xi )+p^0(x),\\
p_1(x,\xi )&=&p_1^1(x,\xi )+p_1^0(x),
\end{eqnarray*}
then we know that $p^2(x,\xi )$ is elliptic; $p^2(x,\xi )\backsim \vert 
\xi \vert ^2$, and 
\begin{eqnarray}\label{heat.24.5.3}
\widetilde{P}(x,\xi ,{1\over \lambda };\widetilde{h})&=&(p^2(x,\xi )+{1\over 
\lambda }p^1(x,\xi )+{1\over \lambda ^2}p^0(x))\\
&&+{h\over \lambda 
}(p_1^1(x,\xi )+{1\over \lambda }p_1(x))+({h\over \lambda 
})^2p_2(x).\nonumber
\end{eqnarray}
\par If $\lambda $ is \sufly{} large, then $p^2(x,\xi )$ is dominating in 
the region $\vert \xi \vert \backsim 1$, and we can construct WKB-solutions 
to \no{heat.24.5.1} ${\rm mod\,}{\cal O}(\widetilde{h}^\infty )$ with all the 
derivatives, of the form 
\ekv{heat.24.5.4}
{
e^{{i}\widetilde{\psi }(s,x,\widetilde{\eta },{1\over 
\lambda })/\widetilde{h}}\widetilde{a}(s,x,\widetilde{\eta },{1\over 
\lambda };\widetilde{h}),
}
with 
$$
\widetilde{\psi }_{\vert s=0}=x\cdot \widetilde{\eta },\ \vert 
\widetilde{\eta }\vert \backsim 1,\ \widetilde{a}_{\vert s=0}=1.
$$
We are now in the elliptic region and it suffices to solve the eikonal \e{} 
and the transport \e{}s to infinte order at $s=0$, since $\Im 
\widetilde{\psi }\backsim s$.

\par If $\eta =\lambda \widetilde{\eta }$, $\vert \widetilde{\eta 
}\vert \backsim 1$, then, at least formally, \no{heat.24.5.4} is just the WKB 
solution $e^{{i\over h}\psi (t,x,\eta )}a(t,x,\eta ;h)$ of the original 
problem \no{heat.2} with ${\psi _\vert}_{t=0}=x\cdot \eta $, 
${a_\vert}_{t=0}=1$, so we can choose 
\begin{eqnarray*}
\psi (t,x,\eta )&=&\lambda \widetilde{\psi }(\lambda t,x,{\eta \over 
\lambda },{1\over \lambda })\\
a(t,x,\eta ;h)&=&\widetilde{a}(\lambda t,x,{\eta \over 
\lambda },{1\over \lambda };{h\over \lambda }),
\end{eqnarray*}
where $\lambda \backsim \vert \eta \vert $. Now $\Im \lambda 
\widetilde{\psi }(\lambda t,x,{\eta \over \lambda },{1\over \lambda 
})\backsim \lambda ^2t$ for $0\le \lambda t\ll 1$ and we get $$ e^{{i\over
h}\psi (t,x,\eta )}={\cal O}(({h\over \lambda })^\infty ),\hbox{ when } \lambda
t\ge (h/\lambda )^{1-\delta }, $$ for any fixed $\delta >0$.

The above discussion indicates how to take care of the uninteresting 
elliptic region. A more complete (and more tedious) treatment could be 
given for example by combining the above scaling argument with a dyadic 
decomposition in $\xi $-space. We observe that $a$ satisfies the symbol 
estimates
$$
\partial _t^k\partial _x^\alpha \partial _\eta ^\beta a={\cal O}(\langle 
\eta \rangle ^{k-\vert \beta \vert }).
$$

\begin{Prop}\label{heat1.1}
\par Modulo a standard reduction to homogeneous non-semi\-classical theory 
(see the 
proof), $U(t)$ is a \fop{} of order 0 with complex phase in the sense of 
\cite{MelSj}, associated to the \ctf{} $\kappa _t$.
\smallskip
\par  $U(t)$ is ${\cal O}(1)e^{-t/C}$ and ${\cal O}_\epsilon 
(1)e^{\epsilon t}$, $\forall \epsilon >0$: $H^s_{\rm 
comp}(\widetilde{X})\to H^s_{\rm loc}(\widetilde{X})$, in the cases $q\ne 
n_-$ and $q=n_-$ respectively.
\smallskip
\par We have 
$$
(h\partial _t+\Delta _q)U(t)={\cal O}(h^\infty )\cases{e^{-t/C},\cr {\cal 
O}_\epsilon (1)e^{\epsilon t},\ \epsilon >0}:\ H^{s-N}_{\rm comp}\to 
H^{s+N}_{loc},$$ in the cases $$\cases{q\ne n_-,\cr q=n_-,}
$$
for all $s\in{\bf R}$, $N\ge 0$.
\end{Prop}

\par\noindent \bf Proof. \rm The statement could be proved directly, but 
it is perhaps more convenient to use the classical theory of \fop{}s with 
complex phase (\cite{MelSj}). The standard trick to get a reduction to 
that 
situation is by adding a variable $x_0$ and to relate \scl{} objects 
(without a tilde) to non-semiclassical objects (with a tilde) in the 
following way: 

\par For \fu{}s we relate the \scl{} ones; $u(x)$, to 
$\widetilde{u}(x_0,x)=e^{ix_0/h}u(x)$.

\par We relate a \scl{} \fop{} 
$$
Fu(x)=\iint e^{{i\over h}\phi (x,y,\theta )}a(x,y,\theta ;h)u(y)dyd\theta 
$$
to a standard (microlocally defined) \fop{}
$$
\widetilde{F}\widetilde{u}(x_0,x)=\iiiint_{\theta _0>0}e^{i(\phi 
(x,y,\theta )\theta _0+(x_0-y_0)\theta _0)}a(x,y,\theta ;{1\over \theta 
_0})\widetilde{u}(y_0,y){dy_0\over 2\pi }dyd\theta _0d\theta ,
$$
so that
$$
\widetilde{F}(e^{{i\over h}x_0}u(x))=e^{{i\over h}x_0}Fu(x).
$$
Here, we require that $\Im \phi \ge 0$, so that the same holds for 
$$
\widetilde{\phi }=\phi (x,y,\theta )\theta _0+(x_0-y_0)\theta _0.
$$
Let $C_\phi =\{ (x,y,\theta );\, \phi '_\theta (x,y,\theta )=0\}$ and 
recall that $\phi $ is \nondeg{} if $d\phi '_{\theta _1}$,...,$d\phi 
'_{\theta 
_N} $ are linearly \indep{}
at every point of $C_\phi $. Then it easy to see that $\phi $ is \nondeg{} 
iff $\widetilde{\phi }$ is, and we have 
$$
C_{\widetilde{\phi }}=\{ (x_0,y_0,\theta _0;x,y,\theta );\, (x,y,\theta 
)\in C_\phi ,\,\, x_0=y_0-\phi (x,y,\theta )\}
$$
We assume (which is the case for $U(t)$) that we are in the \nondeg{} 
case. Then we 
introduce the corresponding canonical relation
$$
\Lambda _\phi =\{ (x,\phi '_x;y,-\phi '_y);\, (x,y,\theta )\in C_\phi \}.
$$
Then for $\widetilde{\phi }$, we have 
\begin{eqnarray*}
\Lambda _{\widetilde{\phi }}&=&\{ (x_0,\xi _0,x,\xi ;y_0,\eta _0,y,\eta 
);\, 
x_0=y_0-\phi (x,y,\theta ),\\ &&\hskip 3truecm \xi _0=\eta _0=\theta _0,\, 
(x,{\xi \over 
\theta _0};y,{\eta \over \theta _0})\in \Lambda _\phi \} .
\end{eqnarray*}
The corresponding relation between the evolution \e{}s is that
$$
(h\partial _t+P(x,hD_x))u={\cal O}(h^\infty )\Leftrightarrow (\partial 
_t+D_{x_0}P(x,D_{x_0}^{-1}D_x))\widetilde{u}=0 \hbox{ \ml{}ly},
$$
when $\widetilde{u}(t,x_0,x)=e^{ix_0/h}u(t,x)$. This is coherent with the 
two other correspondances above, let us just check the geometric one: The 
\ctf{}s associated to $U(t)$, and the solution operator $\widetilde{U}(t)$ 
of the second evolution problem are denoted by $\kappa _t$ and 
$\widetilde{\kappa }_t$
respectively, so that $\widetilde{\kappa }_t$ is obtained by integrating 
the system:
\ekv{heat.24.6}
{
i\dot{x}_0=\partial _{\xi _0}\widetilde{p},\ i\dot{\xi} _0=-\partial 
_{x_0}\widetilde{p},\ 
i\dot{x}=\partial _\xi \widetilde{p},\ i\dot{\xi }=-\partial 
_x\widetilde{p}, 
}
with 
$$
\widetilde{p}(x_0,x;\xi _0,\xi )=\xi _0p(x,\xi _0^{-1}\xi ),
$$
while the corresponding evolution \pb{} giving $\kappa _t$ is
\ekv{heat.24.7}
{
i\dot{x}=\partial _\xi p,\ i\dot{\xi }=-\partial _xp.
} 
Now \no{heat.24.6} becomes 
\begin{eqnarray*}
&&i\dot{x}_0=p(x,\xi /\xi _0)-p'_\xi (x,\xi /\xi _0)\cdot (\xi /\xi _0), 
\ i\dot{\xi }_0=0,\\
&&i\dot{x}=p'_\xi (x,\xi /\xi _0),\ i\dot{\xi}/\xi _0=-p'_x(x,\xi /\xi _0),
\end{eqnarray*}
which reduces to $\kappa _t$
after restriction to $\xi _0=1$.

\par To get the second statement, we observe that $\psi (t,\xi ,\eta 
)\to \psi (\infty ,x,\eta )$ and that the corresponding canonical 
relation $\kappa _\infty $ is strictly positive  with real part being the 
identity relation on $\Sigma $. The statement then follows by the 
description of our operators after conjugation by an FBI-Bargmann \tf{} as 
in \cite{MelSj3}.

\par The proof of the third statement is straight forward. 
\hfill{$\#$}\medskip

\par In the remainder of this section, we assume that $q\ne n_-$

\begin{Prop}\label{heat1.2}
We have
\ekv{heat.24.8}
{[\Delta _q,U(t)]={\cal O}(h^N)e^{-t/C},\ H_{{\rm comp}}^{s-N}\to H_{{\rm 
loc}}^{s+N},}
for all $s\in {\bf R}$ and all $N\ge 0$. 
\end{Prop}

\par\noindent \bf Proof. \rm Using the theory of \cite{MelSj}, we see that 
\ekv{heat.24.9}
{
[\Delta _q,U(t)]u(x)={1\over (2\pi h)^{2n}}\iint e^{{i\over h}(\psi 
(t,x,\eta )-y\cdot \eta )}b(t,x,\eta ;h) u(y)dyd\eta +R(t)u(x),
}
where $R(t)={\cal O}(h^\infty )e^{-t/C}:H_{{\rm comp}}^{-\infty } \to 
H_{\rm 
loc}^\infty $, and $b(t,x,\eta ;h)\sim \sum_0^\infty b_j(t,x,\eta )h^j$ 
satisfies \no{heat.24.2}, \no{heat.24.3} in a region with $\eta $ \bdd{} 
and $\partial _t^k\partial _x^\alpha \partial _\eta ^\beta b={\cal 
O}(\langle \eta \rangle ^{2+k-\vert \beta \vert })$ in a region where $t$ 
is \bdd{}.  Further, we have 
\ekv{heat.24.10}
{
[\Delta _q,U(0)]=0,
}
\ekv{heat.24.11}
{
(h\partial _t+\Delta _q)[\Delta _q,U(t)]={\cal O}(h^\infty )e^{-t/C}.
}
From \no{heat.24.10} we conclude that $b_j(0,\eta )=0$ and from 
\no{heat.24.11} we see that $b_j$
satisfy the same transport \e{}s as $a_j$, and hence 
$$
b_j={\cal O}(t^\infty ),\ b_j={\cal O}(e^{-t/C}{\rm dist\,}(\cdot 
,\Sigma )^\infty),
$$
where we restrict the attention to a region with $\eta $ \bdd{} for 
simplicity.
From this we deduce \no{heat.24.8}.\hfill{$\#$}\medskip

\par Combining the last two propositions, we get 
\ekv{heat.24.12}
{
h\partial _tU(t)+U(t)\Delta _q={\cal O}(h^\infty )e^{-t/C}:H_{\rm comp}
^{-\infty }(\widetilde{X})\to H_{{\rm loc}}^{\infty }(\widetilde{X}).}
From this we get a two-sided parametrix for $\Delta _q$:
\begin{Theo}\label{heat1.3}
We recall that we work with the assumption $q\ne n_-$. Put 
\ekv{heat.24.13}
{
E={1\over h}\int_0^\infty  U(t)dt.
}
Then
\ekv{heat.24.14}
{E={\cal O}(h^{-1}):H_{\rm comp}^s\to H_{\rm loc}^s,}
for every $s\in{\bf R}$, and 
\ekv{heat.24.15}
{
\Delta _qE-1,\, E\Delta _q-1\, = {\cal O}(h^\infty ):H_{\rm 
comp}^{-\infty }\to H_{\rm loc}^\infty .
}
\end{Theo}

\par\noindent \bf Proof. \rm The first estimate follows from the second 
statement in Proposition \ref{heat1.1}. Further,
$$
\Delta _qE={1\over h}\int_0^\infty -h\partial _tU(t)dt+{1\over 
h}\int_0^\infty (h\partial _t+\Delta _q)U(t) dt.
$$
Here the first integral is equal to 
$1$, since $U(0)=1$, and the second integral is ${\cal O}(h^\infty 
):H_{\rm comp}^{-\infty }\to H_{\rm loc}^\infty $ by the last part of 
Proposition \ref{heat1.1}. The proof of \no{heat.24.15} is similar except 
that we use \no{heat.24.12} instead.\hfill{$\#$}\medskip

\section{$\Pi $ as a local projection on ${\cal N}(\Delta _q)$ 
mod$\,{\cal O}(h^\infty )$.}\label{pi}
\par In this section we continue to work in a connected open subset 
where the curvature $\overline{\partial }\partial \phi $ is \nondeg{} of 
signature $(n_+,n_-)$ and we restrict the attention to $(0,q)$-forms, 
with $q=n_-$.

\par Recall that $U(t)$, defined by 
\no{heat.3}, 
is well-defined mod ${\cal O}(h^\infty )$ as an \op{}: $H_{\rm 
comp}^{s}\to H_{\rm loc}^s$ for $t\ge 0$ and as an \op{}: $H_{\rm 
comp}^{s-N}\to H_{\rm loc}^{s+N}$, for $t\ge t_0$ for all $t_0>0$. Put 
\ekv{heat.25}
{
\Pi u={1\over (2\pi h)^{2n}}\iint e^{{i\over h}(\psi (\infty ,x,\eta 
)-y\cdot 
\eta )}a(\infty ,x,\eta ;h)u(y)dyd\eta ,
}
so that $\Pi $
is well-defined ${\rm mod\,}{\cal O}(h^\infty )$ as an \op{} $H_{\rm 
comp}^{s-N}\to H_{\rm loc}^{s+N}$, for all $s\in{\bf R}$, $N\ge 0$. Then 
by 
Proposition \ref{heat1}
, and \no{heat.12}, we have 
\ekv{heat.26}
{
U(t)=\Pi +V(t),\ V(t)={\cal O}(e^{-t/C}):\, H_{\rm comp}^{s-N}\to H_{\rm 
loc}^{s+N},\ t\ge t_0,\ \forall \, t_0>0. 
}
To see this one can introduce $U_s(t)$ with phase $(1-s)\psi (t,x,\eta 
)+s\psi (\infty ,x,\eta )$ and amplitude $(1-s)a(t,x,\eta )+sa(\infty 
,x,\eta ;h)$, $0\le s\le 1$, and show that $\partial _sU_s(t)$ satisfies 
the estimate in \no{heat.26}.
\begin{Prop} \label{heat2} \it We have in the sense of \op{}s: $H_{\rm 
comp}^{s-N}\to H_{\rm loc}^{s+N}$, 
\ekv{heat.27}
{
\Delta _q\Pi \equiv\Pi \Delta _q\equiv 0\ {\rm mod\,}{\cal O}(h^\infty ),
}
\ekv{heat.27.5}
{
\Pi ^*-\Pi \equiv 0\ {\rm mod\,}{\cal O}(h^\infty ),
}
\ekv{heat.28}
{
[\Pi ,V(t)]={\cal O}(e^{-t/C}h^\infty ).
}
\end{Prop}
\par\noindent \bf Proof. \rm We know from \no{heat.21.5}, \no{heat.15} 
that
$$
Z_\phi (e^{{i\over h}\psi (t,x,\eta )}a(t,x,\eta ;h))=e^{{i\over h}\psi 
(t,x,\eta )}b(t,x,\eta ;h),
$$ 
where $b={\cal O}(e^{-t/C})$. Write this as 
$$
{\cal O}(e^{-t/C})=b=e^{-{i\over h}\psi (t,x,\eta )}\circ Z_\phi \circ 
e^{{i\over h}\psi (t,x,\eta )}(a(t,x,\eta ;h))=:Z_{\phi ,\psi ,\eta 
}a(t,x,\eta ;h).
$$
Here 
$$
Z_{\phi ,\psi ,\eta }=Z_{\phi ,\psi _\infty ,\eta } +{\cal O}(e^{-t/C}),
$$
in the sense of Taylor expansions of the coefficients in the $h$-\asy{} 
expansions at $\Sigma $,
by \no{heat.12}. We conclude that 
$$
Z_{\phi ,\psi _\infty ,\eta }a_\infty ={\cal O}(e^{-t/C})+{\cal 
O}(h^\infty ), 
$$
but here the \lhs{} is \indep{} of $t$ and hence
\ekv{heat.29}
{
Z_{\phi ,\psi _\infty ,\eta }a_\infty ={\cal O}(h^\infty ).
}
This means that
\ekv{heat.30}
{
Z_\phi (e^{{i\over h}\psi _\infty }a_\infty )={\cal O}(h^\infty ).
}
Similarly,
\ekv{heat.31}
{
Z_\phi ^*(e^{{i\over h}\psi_\infty  }a_\infty )={\cal O}(h^\infty ).
}
Hence, 
\ekv{heat.32}
{
\Delta _q(e^{{i\over h}\psi_\infty  }a_\infty )={\cal O}(h^\infty ),
}
which implies that 
\ekv{heat.33}
{\Delta _q\Pi ={\cal O}(h^\infty ):\, H_{{\rm comp}}^{s-N}\to H_{\rm 
loc}^{s+N}.}

\par $[\Delta _q,U(t)]$ is an $h$-\fop{} of the same type as $U(t)$. We 
have 
$$
\partial _t[\Delta _q,U(t)]+\Delta _q[\Delta _q,U(t)]=[\Delta _q,\partial 
_tU(t)+\Delta _qU(t)]={\cal O}(h^\infty ) $$ in the sense of such \op{}s
and using also that $[\Delta _q,U(0)]=0$, we get $[\Delta _q,U(t)]={\cal
O}(h^\infty )$ in the sense of such \op{}s and hence $[\Delta
_q,U(t)]={\cal O}(h^\infty )$ as an \op{}: $H_{{\rm comp}}^{s+2}\to H_{\rm
loc}^s$ for $t\ge 0$ and $H_{{\rm comp}}^{s-N}\to H_{\rm loc}^{s+N}$, for
$t\ge t_0>0$.  It follows that $[\Pi ,\Delta _q]={\cal O}(h^\infty ):H_{\rm
comp}^{s-N}\to H_{\rm loc}^{s+N}$ and together with \no{heat.33}, this
gives \no{heat.27}.

\par Next we see that in the sense of $h$-\fop{}s: 
$$
\partial _tU^*+\Delta _qU^*\equiv (\partial _tU+U\Delta _q)^*\equiv 
(\partial _tU+\Delta U)^*\equiv 0.
$$
Hence $\partial _t(U^*-U)+\Delta _q(U^*-U)\equiv 0$, $(U^*-U)(0)=0$, so 
by considering again the transport \e{}s, we get $U^*\equiv U$. It follows 
that $\Pi ^*\equiv \Pi $, so we have \no{heat.27.5}.

\par Consider $[U(t),U(s)]$, $0\le s<\infty $ (after introducing a cutoff 
near the diagonal to make our operators properly supported without 
affecting any other properties). This commutator is 
obviously a \fop{} associated to $\kappa _{t+s}$. For $s=0$, we have 
$$
[U(0),U(t)]=[1,U(t)]=0. 
$$
Moreover, since $\Delta _q$ commutes with $U(t)$:
$$
(h\partial _s+\Delta _q)[U(s),U(t)]\equiv [(h\partial _s+\Delta 
_q)U(s),U(t)]\equiv [0,U(t)]=0.
$$
From considering the transport \e{}s for the amplitude of $[U(s),U(t)]$ 
with the phase $\psi (t+s,x,\eta )-y\cdot \eta $, we see that 
$[U(s),U(t)]=0$. Letting $s\to \infty $, we get $[\Pi ,U(t)]\equiv 0$ and 
\no{heat.28}
follows.
\hfill{$\#$}\bigskip

\par 
For $\Re z<0$, we put
\ekv{pi.1}
{R(hz)=-{1\over h}\int_0^\infty  e^{tz} U(t)dt={\cal O}(h^{-1}):H_{\rm 
comp}^s\to H_{\rm loc}^s.}
Then modulo ${\cal O}(h^\infty ):H_{\rm comp}^{s-N}\to H_{\rm loc}^{s+N}$ 
we have,
\begin{eqnarray*}
\Delta _qR(hz)&=&-{1\over h}\int_0^\infty  e^{tz}\Delta _qU(t)dt\equiv 
{1\over h}\int_0^\infty e^{tz}h\partial _tU(t)dt\\
&=&-\int_0^\infty \partial _t(e^{tz})U(t)dt-1=hzR(hz)-1. 
\end{eqnarray*}
We also have $R(hz)\Delta _q\equiv \Delta _qR(hz)$, so we get
\ekv{pi.2}
{
(hz-\Delta_ q)R(hz)\equiv R(hz)(hz-\Delta _q)\equiv 1.
}
In order to extend to a domain, $\Re z<1/(2C)$, we first rewrite \no{pi.1} 
as 
$$
R(hz)=-{1\over h}\int_0^\infty e^{tz }(\Pi +V(t))dt={1\over hz}\Pi 
-{1\over 
h}\int_0^\infty  e^{tz}V(t)dt,
$$
and for 
\ekv{pi.2.5}
{\Re z<1/(2C),\ \vert z\vert \ge h^{N_0},} with $N_0>0$ \ably{} large but 
fixed, we define
\ekv{pi.3}
{
R(hz)={1\over hz}\Pi -{1\over h}\int_0^\infty e^{tz}V(t)dt={\cal O}(\vert 
hz\vert ^{-1}+h^{-1}):\,H_{\rm comp}^s\to H_{\rm loc}^s.
}
Then this is a \hol{} extension of $R(hz)$, defined by \no{pi.1}. It is 
therefore no surprise that \no{pi.2} remains valid (even though we cannot 
appeal to unique \hol{} extension, since we work with errors that are 
${\cal O}(h^\infty )$): Use that $(h\partial t+\Delta _q)V(t)={\cal 
O}(h^\infty  e^{-t/C}):$$H_{\rm comp}^{s-N}\to H_{\rm loc}^{s+N}$ (cf 
\no{heat.26}, \no{heat.27}) 
$V(0)=1-\Pi $, to get ${\rm mod\,}{\cal O}(h^\infty ): H_{\rm 
comp}^{s-N}\to H_{\rm loc}^{s+N}$,
\begin{eqnarray*}
\Delta _qR(hz)&\equiv & -{1\over h}\int_0^\infty \Delta 
_qV(t)dt\equiv\int_0^\infty e^{tz}\partial _tV(t)dt\\
&=&\Pi -1 -z\int_0^\infty  e^{tz}V(t)dt=\Pi -1-hz{1\over h}\int_0^\infty  
e^{tz}V(t)dt\\
&=& \Pi -1+hz(R(hz)-{1\over hz}\Pi )=hzR(hz)-1,
\end{eqnarray*}
so indeed we have \no{pi.2} for $z$ in the region \no{pi.2.5}.
\begin{Prop} \label{pi1} We have 
\ekv{pi.4}
{
\Pi ={1\over 2\pi i}\int_{\vert z\vert =r}R(z)dz,
}
if $h^{N_0}\le r\le {1/(2C)}$. Moreover,
\ekv{pi.5}
{
\Pi ^2\equiv \Pi .
}
\end{Prop}

\par In order for \no{pi.5} to make sense, we have multiplied the
distribution kernel of $U(t)$ by a cutoff near
the diagonal in order to make all the \op{}s properly supported without 
changing any of their other properties. \medskip

\par\noindent \bf Proof. \rm \no{pi.4} is immediate from \no{pi.3}, since 
the last 
term in \no{pi.3} is \hol{} in $\vert z\vert <1/(2C)$. To prove \no{pi.5}, 
we 
follow the standard procedure and establish first an approximate version 
of the resolvent identity when $h^{N_0}\le \vert z\vert ,\,\vert w\vert 
\le 
h/(2C)$, modulo ${\cal O}(h^\infty )$: $H_{\rm comp}^{s-N}\to H_{\rm 
loc}^{s+N}$,
\ekv{pi.6}
{
R(z)-R(w)\equiv R(z)(w-z)R(w)\equiv R(w)(w-z)R(z).
}
Write 
$$
(z-\Delta _q)-(w-\Delta _q)=(z-w),
$$
and apply $R(z)R(w)$. Then \no{pi.6} follows.

\par Using \no{pi.6}, we write
\begin{eqnarray*}
\Pi ^2&=&({1\over 2\pi i})^2\int_{\vert z\vert =r_1}\int_{\vert w\vert 
=r_2}R(z)R(w)dwdz\\&\equiv
&
({1\over 2\pi i})^2\int_{\vert z\vert =r_1}\int_{\vert w\vert 
=r_2} (w-z)^{-1}R(z)dwdz\\ &&+
({1\over 2\pi i})^2\int_{\vert w\vert =r_2}\int_{\vert z\vert 
=r_1} (z-w)^{-1}R(w)dzdw.
\end{eqnarray*}
Choose $h^{N_0}\le r_1<r_2\le h/(2C)$. In the second integral, we first 
integrate \wrt{} $z$ and get $0$. In the first integral, we first 
integrate 
in $w$ and get 
$$
{1\over 2\pi i}\int_{\vert z\vert =r_1}R(z)dz=\Pi .
$$
\par\noindent \hfill{$\#$}\medskip

\par The next result together with \no{heat.27}, \no{heat.28},
\no{pi.5} 
says that in an approximate sense $\Pi $
is the \og{} projection onto the kernel of $P$ and that $1-\Pi $
is approximately the \og{} projection onto the range of $P$:
\begin{Theo}\label{pi2} For $h^{N_0}\le 
r\le h/(2C)$, put
\ekv{pi.7}
{
E=-{1\over 2\pi i}\int_{\vert z\vert =r}{1\over z}R(z)dz={\cal O}({1\over 
h}):H_{\rm comp}^s\to H_{\rm loc}^s.
}
Then modulo ${\cal O}(h^\infty ):\,H_{\rm comp}^{s-N}\to H_{\rm 
loc}^{s+N}$,
\ekv{pi.8}
{
1\equiv \Pi +\Delta _qE\equiv \Pi +E\Delta _q.
}
\end{Theo}

\par\noindent \bf Proof. \rm Since $E\Delta _q\equiv \Delta _qE$, we only 
have to prove 
the first relation in \no{pi.8}:
\begin{eqnarray*}
\Delta _qE&=&-{1\over 2\pi i}\int_{\vert z\vert =r}{1\over z}\Delta 
_qR(z)dz\\&=&-{1\over 2\pi i}\int_{\vert z\vert =r}{1\over z}(\Delta 
_q-z)R(z)dz-{1\over 2\pi i}\int_{\vert z\vert =r}R(z)dz\\
&\equiv &{1\over 2\pi i}\int_{\vert z\vert =r}{1\over z}dz-\Pi =1-\Pi .
\end{eqnarray*}
\par\noindent \hfill{$\#$}
\medskip \par From the discussion around \no{heat.7}--\no{heat.14}, we 
recollect that \ufly{} for $t\ge t_0>0$:
\ekv{pi.9}
{
\psi (t,x,\eta )=x\cdot \eta +{\cal O}({\rm dist\,}(x,\eta ;\Sigma )^2),
} 
\ekv{pi.10}
{
\Im \psi (t,x,\eta )\sim {\rm dist\,}(x,\eta ;\Sigma )^2.
}
The complex stationary phase method (\cite{MelSj}) then permits us to 
carry out the $\eta $-integration in \no{heat.24.1}, \no{heat.25}, to get 
\begin{Theo}\label{pi3}
For every $t_0>0$, we have \ufly{} for $t\ge t_0$
\ekv{pi.11}
{
U(t)u(x)=h^{-n}\int e^{{i\over h}\widetilde{\psi 
}(t,x,y)}b(t,x,y;h)u(y)m(dy)+R(t)u(x),
}
\ekv{pi.12}
{
\Pi u(x)=h^{-n}\int e^{{i\over h}\widetilde{\psi 
}(\infty ,x,y)}b(\infty ,x,y;h)u(y)m(dy)+R(\infty )u(x),
}
where
\ekv{pi.13}
{
b(t,x,y;h)\sim \sum_0^\infty b_j(t,x,y;h)h^j,
}
\ekv{pi.14}
{
R(t)u(x)=\int r(t,x,y;h)u(y)m(dy),
}
\ekv{pi.15}
{
\partial _{(t,x,y)}^\alpha r={\cal O}(h^\infty ),
}
\ekv{pi.16}
{
\Im \widetilde{\psi }(t,x,y)\sim \vert x-y\vert ^2,\ \widetilde{\psi 
}(t,y,x)=-\overline{\widetilde{\psi }}(t,x,y),
}
\ekv{pi.17}
{
{\rm graph\,}\kappa _t=\{ (x,\partial _x\widetilde{\psi 
}(t,x,y);y,-\partial _y\widetilde{\psi }(t,x,y));\, (x,y)\in{\rm 
neigh\,}({\rm diag\,}(\widetilde{X}\times \widetilde{X}))\} ,
}
\ekv{pi.18}
{
{\partial _{t,x,y}^\alpha (\widetilde{\psi }(t,x,y)-\widetilde{\psi 
}(\infty ,x,y))_\vert }_{y=x}={\cal O}(e^{-t/C})
}
and similarly for $b_j$.
\end{Theo}

\section{The global null-projection}\label{gl}

We first recollect what we have done locally. Let $s$ be a local 
non-vanishing \hol{} section of $L$, defined on $\widetilde{X}\subset X$. 
Write $\vert s(x)\vert ^2=e^{-2\phi (x)}$, and recall that we have the 
unitary map
\begin{eqnarray}\label{gl.0}
{\cal E}^{0,q}(\widetilde{X })\ni u&\mapsto & \widetilde{u}=(se^\phi 
)^ku\in {\cal E}^{0,q}(\widetilde{X};L^k)\\
\overline{\partial }_s&\mapsto & h\overline{\partial }\nonumber\\
\Delta _q&\mapsto & \widetilde{\Delta }_q,\nonumber 
\end{eqnarray}
where $\widetilde{\Delta }_q=h\overline{\partial }h\overline{\partial 
}^*+h\overline{\partial }^*h\overline{\partial }$ is the Hodge Laplacian 
on ${\cal E}^{0,q}(\widetilde{X};L^k)$. Assume the curvature is 
\nondeg{} with $n_-=q$ on $X$. In Section \ref{pi} we constructed an 
approximate resolvent for $\Delta _q$ for $z$ in the domain \no{pi.2.5} 
and 
an approximate null-projection of the form
\ekv{gl.1}
{\Pi u(x)=h^{-n}\int e^{\psi (x,y)/h}b(x,y;h) u(y)m(dy),\ h=1/k,
}
where our new $\psi $ is related to $\widetilde{\psi }(\infty ,x,y)$ in 
\no{pi.12} by 
\ekv{gl.2}
{
\psi (x,y)=i\widetilde{\psi }(\infty ,x,y),
}
so that \no{pi.16}, \no{pi.17} give
\ekv{gl.3}
{
\Re \psi (x,y)\sim -\vert x-y\vert ^2,\ \psi (y,x)=\overline{\psi (x,y)}
}
\ekv{gl.4}
{
d_x{1\over i}\psi (x,y)\in J_+,\ -d_y{1\over i}\psi (x,y)\in J_-.
}
When $x=y$, this implies that 
$$d_xi^{-1}\psi (x,y)=-d_yi^{-1}\psi (x,y)\in \Sigma .$$
We also know from the construction 
that 
\ekv{gl.5}
{
\psi (x,x)=0.
}
On the other hand, we know that $\Sigma $ is given by $\Re \xi dx=\Re 
{2\over i}{\partial \phi \over \partial x}dx$ (using the notations of
Section \ref{bas} but writing $x,\xi$ instead of $z,\zeta$), so we get for 
$x=y$:
$$
d_x{1\over i}\psi (x,y)=\Re {2\over i}{\partial \phi \over \partial 
x}dx={1\over i}{\partial \phi \over \partial x}dx-{1\over i}{\partial \phi 
\over \partial \overline{x}}d\overline{x}.
$$
Hence for $x=y$:
\ekv{gl.6}
{
{\partial \psi \over \partial x}={\partial \phi \over \partial x},\ 
{\partial \psi \over \partial \overline{x}}=-{\partial \phi \over 
\partial \overline{x}},\ 
{\partial \psi \over \partial y}=-{\partial \phi \over \partial x},\ 
{\partial \psi \over \partial \overline{y}}={\partial \phi \over
  \partial \overline{x}}. 
}
Since $\Pi $ is \sa{} modulo ${\cal O}(h^\infty )$, we also have 
\ekv{gl.7}
{
b(x,y;h)^*=b(y,x;h),
} 
where the * indicates that we take the complex adjoint of 
$$
b(x,y;h):\, {\Lambda }^{0,q}T_y^*X\to {\Lambda }^{0,q}T_x^*X.  
$$

\par In terms of 
$$\widetilde{u}=(se^\phi )^ku,\,\widetilde{v}=(se^\phi )^kv\in {\cal 
E}^{0,q}(\widetilde{X};L^k),$$
we get from $v=\Pi u$, that $\widetilde{v}=\widetilde{\Pi }\widetilde{u}$, 
with 
\ekv{gl.7.5}
{
\widetilde{v}=h^{-n}\int e^{\psi (x,y)/h}\widetilde{b}(x,y;h) 
\widetilde{u}(y) m(dy),
}
where the "symbol" 
\ekv{gl.8}
{
\widetilde{b}(x,y;h)=(s(x)e^{\phi (x)})^kb(x,y;h)(s(y)e^{\phi (y)})^{-k}
}
maps
\ekv{gl.9}
{
L_y^k\otimes {\Lambda }^{0,q}T_y^*X \to L_x^k\otimes {\Lambda }^{0,q}T_x^*X
}
and satisfies \no{gl.7}, now in the sense of maps as in \no{gl.9}. Notice 
that though $se^\phi $ is normalized, the "symbol" $\widetilde{b}$ may 
contain oscillations, contrary to the true symbol $b$.

\par Let $s_1$ be a second non-vanishing local \hol{} section of $L$ with 
$\vert s_1\vert =e^{-2\phi _1}$, so that $s_1e^{\phi _1}$ is normalized. 
In the intersection of the domains of definition, we have 
$$
s_1e^{\phi _1}=se^\phi  e^{ig},
$$
with $g$ real and $\phi _1-\phi $ \plh{}. We then have the local 
representation 
$$
\widetilde{u}=(s_1e^{\phi _1})^ku_1,\ \widetilde{v}=(s_1e^{\phi _1})^kv_1 ,
$$
and a null-projection that is unitarily equivalent to the one in 
$\no{gl.1}$:
\ekv{gl.10}
{\Pi_1 u_1(x)=h^{-n}\int e^{\psi_1 (x,y)/h}b_1(x,y;h) u_1(y)m(dy),\ h=1/k.}
Since the heat parametrix constructed in Section \ref{heat} is unique 
${\rm mod\,}{\cal O}(h^\infty ):H_{\rm comp}^{-\infty }\to H_{\rm 
loc}^\infty $, we have the corresponding facts for the local resolvents 
and 
null-projections, so \no{gl.10} necessarily leads to the same relation 
\no{gl.7.5}, and we can also relate $\Pi$, $\Pi_1$ more directly, by 
writing 
$$
u_1=\big({se^\phi \over s_1e^{\phi _1}}\big)^ku=e^{-ikg}u=e^{-ig/h}u,
$$
to get
$$
\Pi _1=e^{-ig/h}\circ \Pi \circ e^{ig/h},
$$
so 
\ekv{gl.11}
{
b_1=b,\ \psi _1(x,y)=\psi (x,y)-ig(x)+ig(y).
}
In particular, $\Re \psi (x,y)$ does not depend on the choice of local 
\hol{} section $s$. The argument above gives a clear idea about 
the \asy{} behaviour of the kernel of the projection onto the space of 
$q$-harmonic forms. To justify this idea we shall consider 
the global resolvents. 

\par On the full \mfld{} $X$ we know that the Hodge Laplacians 
$\widetilde{\Delta }_{q-1}$, $\widetilde{\Delta }_{q+1}$ have no spectrum 
below $h/C$ for some $C>0$ (as could easily be proved using Theorem 
\ref{heat1.3}) and by a standard argument, we conclude that the spectrum 
of $\widetilde{\Delta }_q$ below $h/C$ is reduced to $\{ 0\}$.  For $z$ in
a set \no{pi.2.5} we can glue together the local \op{}s $R(z)$ of Section
\ref{pi} to an \op{} $\widetilde{R}(z)$ (or rather we first glue together 
the locally unique heat kernels to a global one and then define 
$\widetilde{R}(z)$ as in \no{pi.3}) in such a way that 
\ekv{gl.12}
{
(z-\widetilde{\Delta }_q)\widetilde{R}(z)\equiv 
\widetilde{R}(z)(z-\widetilde{\Delta }_q)\equiv 1\ {\rm mod\,}{\cal 
O}(h^\infty ):\, H^{-\infty }(X)\to H^\infty (X).
}
Here we define the Sobolev spaces $H^s(X)=H^s(X,L)$ of sections 
of $L^k$ with $h=1/k$ in a straight forward way from the local 
representations \no{gl.0}, by means of coverings and partitions of unity. 
The choice of such coverings and partitions will affect the $H^s$-norm only 
up to an equivalence that is uniform in $k$. 

\par Since $\widetilde{\Delta }_q$ is an elliptic \op{} in the classical 
sense, we know on the other hand that for $z$ in the set \no{pi.2.5},
\ekv{gl.13}
{
(z-\widetilde{\Delta }_q)^{-1}={\cal O}(h^{-N_0-1}):\, H^s\to H^{s+2}
}
for all $s\in{\bf R}$, so combining this with \no{gl.12}, we get 
\ekv{gl.14}
{
(z-\widetilde{\Delta }_q)^{-1}\equiv \widetilde{R}(z)\ 
{\rm mod\,}{\cal O}(h^\infty ): H^{-\infty }\to H^\infty .
}
Notice that the distribution kernel of an \op{} which is ${\cal 
O}(h^\infty ):H^{-\infty }\to H^\infty $ is ${\cal O}(h^\infty )$ together 
with all its derivatives. On the other hand, the approximate global 
projection $\widetilde{\Pi }$ discussed earlier in this section satisfies 
(cf. Proposition \ref{pi1})
\ekv{gl.15}
{
\widetilde{\Pi }\equiv {1\over 2\pi i}\int_{\vert z\vert 
=r}\widetilde{R}(z)dz\ {\rm mod\,}{\cal O}(h^\infty ):H^{-\infty }\to 
H^\infty ,
}
while the true nullspace projection of $\widetilde{\Delta }_q$,
\ekv{gl.18}
{
\Pi _0:1_{\{ 0\} }(\widetilde{\Delta }_q)
}
satisfies 
\ekv{gl.17}
{
\Pi _0={1\over 2\pi i}\int_{\vert z\vert =r}(z-\widetilde{\Delta 
}_q)^{-1}dz.
}

\par Combining \no{gl.14}, \no{gl.15}, \no{gl.17}, \no{gl.7}, we get the 
main result of this work:
\begin{Theo}\label{gl1}
Let $L$ be a Hermitian \hol{} line bundle over a compact complex \mfld{}
$X$ and fix a positive smooth measure $m(dx)$ on $X$, so that the Hodge
Laplacian $\widetilde{\Delta }_q=\widetilde{\Delta
}_{q,k}=\overline{\partial }^*\overline{\partial }+ \overline{\partial }
\overline{\partial }^*$ is well-defined on $(0,q)$-forms with coefficients 
in $L^k$, $k\in{\bf N}$. Assume the curvature of $L$ has constant 
signature $(n_-,n_+)$ with $n_-+n_+=n:={\rm dim\,}X$. Then for $k\gg 1$, the 
null-space of $\widetilde{\Delta }_{q,k}$ is reduced to $0$ when $q\ne 
n_-$.

\par In the case $q=n_-$, let $s$ be a non-vanishing \hol{} section of $L$ 
on the open subset $\widetilde{X}$, so that \no{gl.0} gives a unitary map 
between $(0,q)$-forms on $\widetilde{X}$ and $(0,q)$-forms on 
$\widetilde{X}$ with coefficients in $L^k$. If $\Pi _0$ denotes the \og{} 
projection onto the null-space of $\widetilde{\Delta }_q$, we put 
$\Pi _{0,s}u=(se^\phi )^{-k}\Pi _0(se^\phi )^ku$, $u\in 
L^2(\widetilde{X},\Lambda ^{0,q}T^*\widetilde{X})$. Then the distribution 
kernel of $\Pi _{0,s}$ is of the form
\ekv{gl.19}
{
K_{\Pi _{0,s}}(x,y)=h^{-n}e^{\psi (x,y)/h}\widetilde{b}(x,y;h)+r(x,y;h),\ 
h=1/k,
}
with $\psi $, $\widetilde{b}$ as in \no{gl.7.5}, \no{gl.2} , \no{pi.12}, \no{gl.4}, 
\no{gl.5}, \no{gl.6}, \no{gl.7} and where 
$\partial _{x,y}^\alpha r={\cal O}(h^\infty )$ for all $\alpha $.
\end{Theo}
\begin{rem}\label{gl2}
\rm Theorem \ref{gl1} also holds for the more general situation of
$(0,q)$-forms with values in $L^{k}\otimes E,$ where $E$ is a rank
$r$ holomorphic Hermitian vector bundle over $X$. Indeed, locally
$E$ is isomorphic to the trivial holomorphic vector bundle $C^{r}\times X$
with a Hermitian metric $\gamma.$ The local expression \no{db.20}
for $\Delta_{q}$ then still holds if the operators $hZ_{j}+Z_{j}(\phi)$
and their adjoints are tensored by $I_{r},$ the identity matrix on
$C^{r}.$ This follows from the fact that the Hermitian metric $\gamma$
on $E$ is independent of $h=k^{-1}.$ Moreover, globally there is
still a spectral gap for the same reason (as is well-known), giving
an asymptotic expansion as before. For example, if $\mu_{n}$ is a
general volume form on $X$ and $\omega_{n}$ is the one induced by
the given Hermitian metric on $X,$ then the function $\mu_{n}/\omega_{n}$
defines a Hermitian metric on the trivial line bundle $E.$
\end{rem}

\section{\label{ccs}Change of complex structure}

In this section we will investigate some relations to \cite{s-z}
(see also \cite{BoGu}). Let us first recall the setting in \cite{s-z}.
Assume given a symplectic manifold $(X,\omega)$ such that $\pi\omega$
represents an integral cohomology class. Then there exists a Hermitian
line bundle $L$ over $X$ with a unitary connection $\nabla$ whose
curvature satisfies $\frac{i}{2}\Theta=\omega$ (compare Section 
\ref{sub:connections}
for notation). Take an almost complex structure $J'$ on $X$ (i.e.
$J'\in\textrm{End}(TX),\, J'^{2}=-I)$ such that \begin{equation}
\begin{array}{lrcl}
(i) & \omega(J'v,J'w) & = & \omega(v,w)\\
(ii) & \omega(v,J'v) & > & 0\end{array}\label{ccs.1}\end{equation}
 for all $v,w$ in $TX.$  We decompose\[
TX\otimes\C=T^{1,0}(X,J')\oplus T^{0,1}(X,J'),\]
so that $J'\otimes\C=i\oplus-i$ (then \ref{ccs.1} means
that $\omega$ is a positive $(1,1)-$form with respect to $J').$ Then we 
get an operator 
$\nabla^{0,1}:=\overline{\partial}_{J'}$
acting on sections of $L.$ Furthermore, a Riemannian metric $g$ is
said to be compatible with $J'$ if \begin{equation}
g(J'v,J'v)=g(v,v),\label{ccs.2}\end{equation}
i.e. $g$ corresponds to the real part of a Hermitian metric on 
$T^{1,0}(X,J').$

In
\cite{s-z} Shiffman and Zelditch, motivated by the work $\cite{Don}$
of Donaldson, define a sequence of spaces imitating $H^{0}(X,L^{k})$
in the usual integrable case. A naive choice would be the kernel of
$\overline{\partial}_{J'}$ acting on $L^{k},$ but if $J'$ is non-integrable
then these spaces are too small. Instead, Shiffman and Zelditch, following
Boutet de Monvel and Guillemin \cite{BoGu}, introduce a sequence
of spaces of so called \emph{asymptotically almost holomorphic sections.}
The main result in the present section (theorem $\ref{PartSerre}$)
says that the dimension of the null-space of $\Delta_{q}$, studied
in the previous sections, coincides with the dimension of a space
of asymptotically almost holomorphic sections. The latter space is
defined with respect to a new almost complex structure on the original
complex manifold $X.$ It would be very interesting to know if this
correspondence could be extented to the level of Bergman kernels in
a suitable sense, in particular in view of the results in
$\cite{MaMar2}$ on
lower order terms of generalized Bergman kernels. It should finally
be pointed out that in \cite{s-z} the analysis is reduced to the
homogenous theory in \cite{BoGu} by adding a varible dual to $k,$
i.e. by embedding $X$ in the unit circle bundle in $L^{*}$ (this
is a global version of the reduction used in proposition $\ref{heat1.1}$).
But since we work directly in a semiclassical inhomogenous setting
we have developped some of the material in \cite{s-z} from our point
of view.

\subsection{\label{sub:The-pair-}The pair $(J,J')$ }

We now return to the situation in the previous chapters, i.e. we take
$L$ to be a Hermitian line bundle
which is also holomorphic over $(X,J)$ where $J$ denotes the integrable
complex structure. Then it has a canonical connection $\nabla$ (see
Section \ref{sub:connections}). The curvature $\Theta$ of $\nabla$
is assumed to be of signature $(n_-,n_+)=(q,n-q)$ and we will
call $q$ the \emph{index} of $\Theta.$ Hence, $\omega:=\frac{i}{2}\Theta$
is not positive with respect to $J,$ unless $q=0.$ However, given
a Hermitian metric $H$ on $T^{1,0}(X,J)$ as in Section \ref{bas}
(so that its real part corresponds to $g$ in (\ref{ccs.2}))
we can define an almost complex structure $J'$ making $\omega$ positive,
in the following way. Split the real tangent bundle $TX$ as 
\begin{equation}
TX=(TX)_{-}\oplus(TX)_{+}\label{ccs.3}\end{equation}
 according to the positive and negative eigenspaces of 
$\omega(\cdot,J\cdot)$
with respect to the metric $g.$ Then $J$ splits as $J_{+}\oplus J_{-}$
by restriction. Now define $J'$ by the splitting \begin{equation}
J'=(-J_{-})\oplus J_{+}.\label{ccs.4}\end{equation}
Then, clearly, $\omega(v,J'v)>0.$ Equivalently, let $e^{i}$ be a
local frame for $T^{*0,1}(X,J),$ orthonormal with respect to $H,$
such that \begin{equation}
\Theta=\sum_i \lambda_{i}\overline{e^{i}}\wedge e^{i},\label{ccs.5}\end{equation}
 where $\lambda_{i}<0$ for $i\leq q$ and $\lambda_{i}>0$ for $i>q.$
Let $e'^{i}=\overline{e^{i}}$ for $i\leq q$ and $e'^{i}=e^{i}$
for $i>q.$ Then $T^{*0,1}(X,J')$ is spanned by all $e'^{i}$ and
$\Theta=\left|\lambda_{i}\right|\overline{e'^{i}}\wedge e'^{i}$ satisfies
(\ref{ccs.1}). The canonical connection $\nabla$ on
the Hermitian line bundle $L$ induced by $J$ now gives an operator
$\nabla^{0,1}:=\overline{\partial}_{J'}$ (decomposing with respect
to $J')$. 

\par In the sequel $X$ and $X'$ will denote the almost complex manifolds
$(X,J)$ and $(X,J'$) respectively and in general a prime on an object
will indicate that it is defined with respect to the almost complex
structure $J'.$ 

\begin{rem}
\label{ccs1}\rm Even though the pair $(\omega,J')$ fits into
the setup of \cite{s-z} it should be pointed out that the Riemannian metric
$\omega(v,J'w)$ on $X$ was used in \cite{s-z}, but we will use
the 
the
Riemannian metric $g$ induced by the given Hermitian metric $H$
instead. It should be pointed out that the results in this paper are
independent of the metric, but the metric may be important in a more
refined study involving Bergman kernels.
Also, in \cite{s-z} the asymptotics
of projection operators acting on $L^{k}$ were studied, but as mentioned
there, the case $L^{k}\otimes E$ where $E$ is a complex vector bundle
is similar. In Section \ref{PartSerre} we will study
$L^{k}\otimes E$ for a certain complex line bundle $E=K_{X'}^{-}.$ 
\end{rem}

\subsection{\label{sub:connections}Connections and commutation relations}

Let us first recall some basic facts about connections \cite{We},\cite{gr}.
A \emph{connection} $\nabla$ on a complex line bundle $L$ over a
real manifold $X$ is an operator \[
\nabla : C^{\infty}(X;L)\rightarrow C^{\infty}(X;L\otimes 
T^{*}X)\]
 satisfying Leibniz rule: $\nabla(fs)=df\otimes s+f\nabla s$ for
$f$ a function and $s$ a section of $L.$ Given a vector field $v$
on $X$ the contracted operator $\nabla_{v}$ on sections of  $L$ is
called the covaraint derivative along $v.$ The curvature two-form
$\Theta$ of $\nabla$ can be defined by\begin{equation}
\Theta(v,w)=[\nabla_{v},\nabla_{w}]-\nabla_{[v,w]},
\label{ccs.6}\end{equation}
where $v$ and $w$ are vector fields on $X.$ If $L$ has a Hermitian
metric $\left\langle \cdot,\cdot\right\rangle ,$ then a connection
$\nabla$ is called \emph{unitary} if \begin{equation}
d\left\langle s,t\right\rangle =\left\langle \nabla s,t\right\rangle 
+\left\langle s,\nabla t\right\rangle \label{ccs.7}\end{equation}
and if $L$ is a holomorphic line bundle over a complex manifold $X,$
then $\nabla$ is called \emph{holomorphic} if \begin{equation}
\nabla^{0,1}=\overline{\partial}\label{ccs.8}\end{equation}
i.e. $\nabla^{0,1}s=0$ if $s$ is a holomorphic section. There is a unique 
unitary holomorphic connection (see below) on a
Hermitian holomorphic line bundle $L.$ If $(X,J)$ is only an almost
complex manifold, any given connection $\nabla$ defines an operator
$\overline{\partial}_{J}:=\nabla^{0,1}$ acting on sections with values
in $L,$ but there is no canonical operator $\overline{\partial}_{J}$
on $L.$ In the following we will only consider unitary connections
$\nabla$ on $L$ over an almost complex manifold $(X,J).$ 

\par The local situation is as follows. Let $t$ be a local unitary trivializing
section of $L$ and let $A$ be the local one form defined by $\nabla 
t=A\otimes t.$
Note that $\nabla$ is unitary (i.e. (\ref{ccs.7}) holds)
precisely when $A$ is an imaginary one form. Now we get the local
representation $\nabla=d+A,$ i.e.\begin{equation}
\nabla(fs)=(d+Af)s,\label{ccs.9}\end{equation}
 and the curvature two-form $\Theta$ of $\nabla$ is locally given
by  \begin{equation}
\Theta=dA.\label{ccs.10}\end{equation}
If $\widehat{t}=e^{ig}t$ is another unitary frame for $L$ over $U,$
then, using Leibniz rule, the corresponding one form is given by 
\begin{equation}
\widehat{A}=A+idg,\label{ccs.11}\end{equation}
 confirming that the curvature two-form (\ref{ccs.10})
is independent of the local frame.
Take local dual orthonormal frames $Z_{i}$ and $e_{i}$ for $T^{0,1}(X,J)$
and $T^{*0,1}(X,J),$ respectively as in Section \ref{bas}. Splitting
$\nabla=\nabla^{1,0}+\nabla^{0,1}$ we may then write \label{*16} \begin{equation}
\nabla=\sum_i (\overline{e^{i}}\nabla_{\overline{i}}+e^{i}\nabla_{i}),
\label{ccs.12}\end{equation}
 where $\nabla_{i}:=\nabla_{Z_{i}}$ and $\nabla _{\overline{i}}$ are the
 corresponding covariant derivatives along $Z_{i}$ and $\overline{Z}_i$ 
 respectively. Let us now consider some
 local commutation relations.  Write \begin{equation}
\begin{array}{rcl}
[Z_{i},\overline{Z_{j}}] & = & 
\sum_p (a_{ij}^{p}Z_{p}-\overline{a_{ji}^{p}}\overline{Z_{p}})\\
{}[Z_{i},Z_{j}] & = & 
\sum_p (f_{ij}^{p}Z_{p}+\overline{N_{ij}^{p}}\overline{Z_{p}})\end{array},
\label{ccs.13}\end{equation}
where the bracket denotes the commutator between the corresponding
differential operators. Then $N_{ij}^{p}$ is identically zero precisely
when $T^{1,0}(X,J)$ is closed under the bracket, which in turn is
equivalent to $J$ being integrable \cite{Ho,s-z}. In general the
$N_{ij}^{p}$ define the so called Nijenhuis tensor of the almost
complex structure $J.$ Now, using formulas (\ref{ccs.6})
and (\ref{ccs.13}) we get the following commutation relations:
\begin{equation}
\begin{array}{rcl}
[\nabla_{i},\nabla_{\overline{j}}] & = & 
\Theta(Z_{i},\overline{Z_{j}})+\sum_p(a_{ij}^{p}\nabla_{p}-
\overline{a_{ji}^{p}}\nabla_{\overline{p}})\\
{}[\nabla_{i},\nabla_{j}] & = &
\sum_p(f_{ij}^{p}\nabla_{p}+\overline{N_{ij}^{p}}\nabla_{\overline{p}}),\end{array}
\label{ccs.14}\end{equation}
 where we have used that $\Theta$ vanishes on $T^{0,1}X\otimes T^{0,1}X,$
by the assumption (\ref{ccs.1}) $(i)$ on $J.$

\par Let us now specialize to our original situation (compare Section 
\ref{sub:The-pair-}),
where $J$ is integrable and $\Theta$ has index $q$ and take a local
frame $e^{i}$ diagonalizing $\Theta.$ Recall (Section \ref{bas})
that $\left|s\right|^{2}=e^{-2\phi}$ where $s$ is a local holomorphic
trivializing section of $L$ so that $t:=:e^{\phi}s$ is a local unitary
section. Then if $\nabla$ denotes a connection satisfying (\ref{ccs.7})
and (\ref{ccs.8}) we see that $\nabla$ is unique since
the local one form $A$ is given by

\begin{equation}
A=-\partial\phi+\overline{\partial}\phi,\,\,\,
\Theta=dA=2\partial\overline{\partial}\phi
\label{ccs.15}\end{equation}
Indeed, the assumption (\ref{ccs.8}) gives as in Section
\ref{bas} that $\nabla^{0,1}$ is locally represented (with respect
to $t$) by\[
\overline{\partial}+\overline{\partial}\phi=\sum_i e^{i}(Z_{i}+Z_{i}\phi)\]
 and since $A$ is imaginary (by $(ii))$ we get (\ref{ccs.15}).
Moreover, we get that $\Theta_{ij}=\lambda_{i}\delta_{ij}$ and 
$N_{ij}^{p}=0$
in (\ref{ccs.14}). Next, introducing the complex structure
$J'$ defined above corresponds to letting \begin{equation}
\nabla
_{i}^{'}=\nabla_{\overline{i}},\, i\leq q\,\textrm{\,\, 
and$\,\,\,\nabla_{i}^{'}=\nabla_{i},\, i>q$}\label{ccs.16} 
\end{equation}
since the decomposition (\ref{ccs.12}) of $\nabla$ changes
in the corresponding way.

\subsection{Symbols and\label{SymbIdeals} ideals}

We will now consider an arbitrary almost complex structure $J$ again
and replace $L$ by $L^{k}$ and consider semiclassical symbols as
in Section \ref{bas} (setting $h=k^{-1}).$ The discussion will
be local on $U$, given a unitary trivializing section $t$ over $U$
and dual orthonormal frames $Z_{i}$ and $e_{i}$ for $T^{0,1}(U,J)$
and $T^{*0,1}(U,J),$ respectively. Any given connection $\nabla$
on $L$ induces a connection on $L^{k}$, that we also denote by $\nabla$
(i.e. locally $A$ in (\ref{ccs.9}) is replaced by $kA).$
We denote by $\sigma$ the semiclassical principal symbol map
(compare the discussion
about semiclassical principal symbols following the proof of proposition
$\ref{db0}$)
and let \[
q_{i}:=\sigma(h\nabla_{i})\]
in terms of the decomposition (\ref{ccs.3}) (i.e. $q_{i}$
is the principal symbol of the $i$th component of 
$h\overline{\partial}_{J}).$
We will call $\mathcal{J}=(q_{1},...,q_{n})$ the \emph{symbol ideal}
of $\overline{\partial}_{J}.$ Since $\nabla$ is unitary, i.e. it
satifies (\ref{ccs.7}), integration by parts gives \begin{equation}
\sigma(-h\nabla_{\overline{i}})=\sigma(h\nabla_{i}^{*})=\overline{q_{i}},
\label{ccs.17}\end{equation}
 also using the general fact that $\sigma(D^{*})=\overline{\sigma(D)}$
in the last equality. Recall the following general relation between
the operator bracket and the Poisson bracket: \[
\sigma[D_{1},D_{2}]=-ih\left\{ \sigma(D_{1}),\sigma(D_{2})\right\} ,\]
 Hence, the commutator relations (\ref{ccs.14}) become:
\begin{equation}
\begin{array}{rcl}
i\left\{ q_{i},\overline{q_{j}}\right\}  & = & 
\Theta(Z_{i},\overline{Z_{j}})+\sum_p(a_{ij}^{p}q_{p}-
\overline{a_{ji}^{p}}\overline{q_{p}})\\
{}i\left\{ q_{i},q_{j}\right\}  & = & \sum_p(
f_{ij}^{p}q_{p}+\overline{N_{ij}^{p}}\overline{q_{p}}),\end{array}
\label{ccs.18}\end{equation}
Let now $J$ and $J'$ be as in Section \ref{sub:The-pair-}. Note
that when $\nabla$ is the canonical connection determined by $J,$
then the local expression (\ref{ccs.15}) shows that \[
h\nabla_{i}=hZ_{i}+hZ_{i}\phi\]
and that (\ref{ccs.18}) is consistent with the formula
(\ref{db.27}). Now (\ref{ccs.18}) and (\ref{ccs.17})
give $q_{i}'=-\overline{q_{i}}$ for $i\leq q$ and $q_{i}'=q_{i}$
for $i>q.$ In particular, the zero varieties in $T^{*}U$ defined by the
symbol ideals $\mathcal{J}$ and $\mathcal{J}'$ coincide and are equal to
the real characteristic variety $\Sigma=\left\{ 
p_{0:}:=\sigma(\Delta_{q})=0\right\} ,$
where $p_{0}$ is as in formula (\ref{db.21}). In Section \ref{heat}
the local almost holomorphic manifold $J^{+}$ in the almost 
complexification
of $T^{*}X$ was introduced. It corresponds to a local ideal 
$\mathcal{J}^{+}$
of local smooth functions on the symplectic manifold $T^{*}X$ such
that for all $f$ in $\mathcal{J}^{+}$\[
\overline{\partial} \widetilde{f}_{J^{+}}(x)=\mathcal{O}(\textrm{Im}x)^{\infty},\]
where $\widetilde{f}$ denotes an almost holomorphic extension of
$f$ from $T^{*}X.$ The properties of $J^{+},$ reviewed in \no{heat.14.5},
when formulated in terms of the ideal $\mathcal{J}^{+}$, can be stated
as the following lemma \cite{BoGu},\cite{s-z},\cite{MeSj1}, where
$\mathcal{I}_{\Sigma}$ denotes the ideal of elements in 
$\mathcal{C^{\infty}}(X,\C )$
vanishing on $\Sigma$.

\begin{lem}
\label{ccs2}There exists a unique positive Poisson
ideal $\mathcal{J}^{+}$ with respect to $\Sigma$ containing $p_{0}.$
That is, there exists a unique ideal 
$\mathcal{J}^{+}\subset\mathcal{I}_{\Sigma}$
with common zero set $\Sigma$ satisfying $(i)$ $\mathcal{J}^{+}$
is closed under the Poisson bracket and $(ii)$ there are generators
$q_{i}$ of $\mathcal{J}^{+}$ such that the matrix $\frac{1}{i}\left\{ 
q_{i}\overline{,q_{j}}\right\} $
is positive definite on $\Sigma$ and $p_{0}\in\mathcal{J}^{+}.$
\end{lem}
Note that by (\ref{ccs.18}) the ideal $\mathcal{J}$
fails to satisfy the positivity condition $(ii)$ above, since $\Theta$
is assumed to have index $q.$ On the other hand, the positive
ideal $\mathcal{J}'$ only satisfies condition $(i)$ mod 
$\mathcal{I}_{\Sigma}$
(compare Proposition \ref{ccs5}). By the uniqueness
of $\mathcal{J}^{+},$ we then deduce that $\mathcal{J}^{+}=\mathcal{J}'$
mod $\mathcal{I}_{\Sigma}^{2}.$  In fact, the ideal $\mathcal{J}^{+}$
can be constructed from $\mathcal{J}'$ by induction with respect
to $N$ on the vanishing order $\mathcal{I}_{\Sigma}^{N}$ \cite{s-z}
so that $\mathcal{J}^{+}$ is unique mod $\mathcal{I}_{\Sigma}^{N}$
for each $N.$

\begin{rem}
\label{ccs3}\rm The uniqueness mod $\mathcal{I}_{\Sigma}^{2}$
in Lemma \ref{ccs2} is equivalent to the well-known
fact that given a Riemannian metric $g$ on a symplectic manifold
$(X,\omega)$ there is a unique almost complex structure $J$ such
that (\ref{ccs.2}) and (\ref{ccs.1}) hold. The
point is that given $(X,\omega,g,J)$ we get a map\[
(T^{0,1}X,J)\rightarrow\mathcal{J}^{+}/\mathcal{I}_{\Sigma}^{2},\,\,\, 
Z_{i}\mapsto\sigma(h\nabla_{Z_{i}})=:q_{i}\]
 By (\ref{ccs.18}) the conditions (\ref{ccs.1})
on $\omega$ correspond to the conditions in Lemma \ref{ccs2}
on the Poisson brackets when restricted to $\Sigma$ and the condition
(\ref{ccs.2}) on $g$ corresponds to the condition on 
$p_{0}=\sigma(\Delta).$
However, Lemma \ref{ccs2} applies to a more general
situation where $p_{0}$ is a general function on a symplectic manifold
$Y$ (replacing $T^{*}X$ with its usual symplectic form) vanishing
to second order on $\Sigma:=\left\{ p_{0}=0\right\} .$ Then one gets
a complex structure on the normal bundle $TY/T\Sigma$ of $\Sigma$
in $Y$ \cite{BoGu}.
\end{rem}
One final
\begin{rem}
\label{ccs4}\rm In Section \ref{heat} $J^{+}$ was only
locally defined, but it corresponds to a global submanifold of the
almost complexification of the affine bundle $AX$ defined in Section
\ref{ax}.
\end{rem}

\subsection{$J'$ is generically non-integrable}

We will call a function $f$ on $T^{*}X$ \emph{fiber affine} if it
is affine on each fixed fiber of $T^{*}X.$ Equivalently, $f$ is fiber
affine if it is the semiclassical principal symbol of a first order
$h-$differential operator on $X.$ Consider $\C^{2}$ with its standard
complex structure and metric and let $L$ be the trivial holomorphic
line bundle with fiber metric $\phi.$ We will also assume that the
index of the curvature $\Theta=2\partial\overline{\partial}\phi$
is one. 

\begin{prop}
\label{ccs5}The almost complex structure $J'$ is
non-integrable for generic fiber metrics $\phi.$ More precisely,
$J'$ is non-integrable if \begin{equation}
\frac{\partial^{3}\phi}{\partial^{2}z_{1}\partial\bar{z}_{2}}\neq0
\label{ccs.19}\end{equation}
 at $0.$  In particular, the ideal $\mathcal{J}^{+}$ has no fiber
affine generators then.
\end{prop}
\par\noindent \bf Proof. \rm
We will identify $\Theta$ with a Hermitian matrix: 
$\Theta_{ij}:=\Theta(\frac{\partial}{\partial\bar{z}_{i}},\frac{\partial}{\partial 
z_{j}})=-2\frac{\partial^{2}\phi}{\partial z_{j}\partial\bar{z}_{i}}$
with respect to the standard orthogonal frame 
$\frac{\partial}{\partial\bar{z}_{i}}$
and we may assume that $\Theta(0)$ is diagonal. Denote by $Z_{i}$
an orthonormal frame diagonalizing $\Theta$ close to $z=0$, i.e.
$D_{ij}:=\Theta(Z_{i},\overline{Z_{j}})=:-\delta_{ij}\lambda_{i}.$
Equivalently, $Z_{i}=U\frac{\partial}{\partial\bar{z}_{i}}$ where
the matrix valued function $U$ satisfies\begin{equation}
\begin{array}{lrcllrcl}
(i) & U^{*}U & = & I & (ii) & U^{*}\Theta U & = & 
D\end{array},\label{eq:u}\end{equation}
 denoting by $U^{*}$ the Hermitian adjoint $\overline{U}^{t}.$ By
the definition (\ref{ccs.4}) of $J'$ and the subsequent discussion
we may take $Z'_{1}=\overline{Z_{1}}$ and $Z'_{2}=Z_{2}.$ In particular,
$J'$ is non-integrable if $a_{21}^{1}$, defined with respect
to $J$ in (\ref{ccs.13}), is non-vanishing at the origin. Now observe
that at $z=0,$ \begin{equation}
-a_{21}^{1}=\left\langle [\overline{Z_{1}},Z_{2}],Z_{1}\right\rangle 
=(\frac{\partial}{\partial z_{1}}u_{21})(0)
\label{ccs.21}\end{equation}
Indeed, $Z_{1}=\frac{\partial}{\partial\bar{z}_{1}}$ at $z=0$ and
when calculating \[
[\overline{Z_{1}},Z_{2}]=[\overline{u_{11}}\frac{\partial}{\partial 
z_{1}}+\overline{u_{12}}\frac{\partial}{\partial 
z_{2}},u_{21}\frac{\partial}{\partial\bar{z}_{1}}+u_{22}
\frac{\partial}{\partial\bar{z}_{2}}],\]
we can use Leibniz rule for the bracket to expand the right hand side
and get terms of the form\[
\overline{(u_{11}}[\frac{\partial}{\partial 
z_{1}},u_{21}])\frac{\partial}{\partial\bar{z}_{1}}+...\]
But since, $u_{ij}(0)=\delta_{ij}$ the other term proportional to
$\frac{\partial}{\partial\bar{z}_{1}}$ vanishes at $z=0,$ proving
(\ref{ccs.21}). Hence, we just have to show that 
$(\frac{\partial}{\partial z_{1}}u_{21})(0)\neq0,$
if (\ref{ccs.19}) holds. To this end, apply 
$\frac{\partial}{\partial z_{1}}$
to (\ref{eq:u}) and use that $U(0)=I$ and $\Theta(0)=D$ to get at
$z=0$\begin{equation}
(i)\ \frac{\partial}{\partial z_{1}}(U^{*}) =  
-\frac{\partial}{\partial z_{1}}U,\quad  (ii)\ \frac{\partial}{\partial 
z_{1}}\Theta+[D,  \frac{\partial}{\partial z_{1}}U]= 
\frac{\partial}{\partial z_{1}}D.\label{eq:uprim}\end{equation}
In particular, \[
\frac{\partial}{\partial 
z_{1}}\Theta_{21}-(\lambda_{1}-\lambda_{2})\frac{\partial}{\partial 
z_{1}}u_{21}=0,\]
i.e. at $z=0$ we have \[
\frac{\partial}{\partial 
z_{1}}u_{21}=2\frac{\partial^{3}\phi}{\partial^{2}z_{1}\partial\bar{z}_{2}}/(\lambda_{2}-\lambda_{1})\]

\par By (\ref{ccs.21}) this proves first part of the proposition 
about the non-integrabili\-ty of $J'.$

\par The second part is a direct consequence of the first part, by the
way $\mathcal{J}^{+}$ is constructed in \cite{BoGu}. Indeed, by the
uniqueness property mod $\mathcal{I}_{\Sigma}^{2}$ in Lemma 
\ref{ccs2}
we have, since $\mathcal{J}^{+}$ is assumed to be generated by fiber
affine functions, that $\mathcal{J}^{+}=(q_{1}^{'},q_{2}^{'},...,q'_{n}).$
But then the assumption that $\mathcal{J}^{+}$ is a Poisson ideal
forces $N_{ij}^{'p}=0$ in the relations corresponding to (\ref{ccs.18})
for the almost complex structure $J'.$ But this contradicts the first
part of the proposition.  \par \hfill{$\#$}

\subsection{Partial Serre duality\label{PartSerre}}

In this section $J$ and $J'$ (and $X$ and $X'$) will be as in
Section \ref{sub:The-pair-}. Moreover, we will consider globally
defined symbol ideals etc (compare Remark \ref{ccs4}).
Denote by $\mathcal{H}^{0,q}(X,L^{k})$ the global null space of 
$\Delta_{q}$
and denote in this section by $\Pi_{X}^{q}$ the orthogonal projection
on $\mathcal{H}^{0,q}(X,L^{k}).$ In simple cases, e.g. when $X$ is a 
product of complex curves
and $L$ is the product of pulled back line bundles 
one can show that that,
for $k$ sufficiently large, any element $\alpha$ in 
$\mathcal{H}^{0,q}(X,L^{k})$
may be written locally as\[
\alpha=fe^{1}\wedge...\wedge e^{q}\]
 with respect to a frame as in formula (\ref{ccs.5}). Moreover
$f$ is holomorphic with respect to a new integrable complex structure
of the form $J'$. 
In fact, this follows from {}``partial
Serre duality'', i.e. Serre duality along the negative directions
of the line bundle in the product case.
In this section we will show that a version of
this phenomenon, with $J'$ possibly non-integrable, persists for general
$X.$

\par Denote by $K_{X}$ the canonical line bundle on $X=(X,J),$ i.e. $K_{X}$
is the holomorphic line bundle $\Lambda^{n,0}(T^{*}X,J)$ (considering
$(T^{*}X,J)$ as a holomorphic vector bundle). The splitting (\ref{ccs.3})
then induces a decomposition of complex line bundles \begin{equation}
K_{X}=K_{X}^{-}\otimes K_{X}^{+}\label{ccs.23}\end{equation}
 where $K_{X}^{-}:=\Lambda^{q,0}((T^{*}X)_{-},J)$ and 
$K_{X}^{+}:=\Lambda^{n-q,0}((T^{*}X)_{+},J).$ Note
that the bundles $K_{X}^{\pm}$ are not holomorphic in general. Now
\[
K_{X'}^{-}:=\overline{K_{X}^{-}}=\Lambda^{q,0}((T^{*}X)_{-},J')\] is a
complex line bundle on $X'$ with a connection induced by the canonical
connection on $(TX,J)$ determined by the metric $g$ and the complex
structure $J.$ Given a sufficiently large integer $k$ the complex line
bundle $L^{k}\otimes K_{X'}^{-}$ over $X'$ has positive curvature with
respect to $J'$ and fits into the setup in the beginning of Section
\ref{ccs}.  Denote by $H^{0}(X',L^{k}\otimes K_{X'}^{-})$ the space of
\emph{asymptotically almost holomorphic sections} defined in \cite{s-z}
(see Remark
\ref{ccs1})
We will just recall that $H^{0}(X',L^{k}\otimes K_{X'}^{-})$ is defined
as the range of a global projection operator $\Pi_{X'}^{0},$ which
is a Fourier integral operator with complex phase and its canonical
relation can be described in the following way. Let $\Sigma'$ be
the real characteristic variety of $\overline{\partial}_{J'}$ and
let $\mathcal{J}'^{+}$ be the ideal obtained from Lemma \ref{ccs2}
applied to $\Sigma'.$ Then the canonical relation may be written
as 
$\mathcal{C}'\mathcal{=J}'^{+}\times_{\Sigma'}\overline{\mathcal{J}'^{+}}$,
which is to be interpreted in terms of bicharacteristic strips as
in the expression for $C_{\infty}$ (the canonical relation of 
$\Pi_{X}^{q})$
in Section \ref{heat}. But since $\Sigma'=\Sigma$, as observed
in Section \ref{SymbIdeals}, the uniqueness in Lemma
\ref{ccs2} gives $\mathcal{J}'^{+}=\mathcal{J}^{+}$
and hence $\mathcal{C}'=\mathcal{C_{\infty}}.$ The construction in
\cite{s-z} actually only determines $\Pi_{X'}^{0}$ mod 
$\mathcal{O}(k^{-\infty})$
i.e. the asymptotics of its distribution kernel is only determined
up to terms of order $\mathcal{O}(k^{-\infty})$. But this means that
the dimension of $H^{0}(X',L^{k}\otimes K_{X'}^{-})$ is independent
of the construction for $k$ sufficiently large. We will now prove
the following  

\begin{thm}
Assume that the index of $\Theta$ is $q$. Then, for \label{thm:serre}$k$
sufficiently large, \begin{equation}
\dim \mathcal{H}^{0,q}(X,L^{k})=\dim H^{0}(X',L^{k}\otimes 
K_{X'}^{-}).\label{ccs.24}\end{equation}
Furthermore, if $K_{X}$ has a square root $K_{X}^{1/2},$ then 
\begin{equation}
\dim\mathcal{H}^{0,q}(X,L^{k}\otimes K_{X}^{1/2})=\dim 
H^{0}(X',L^{k}\otimes K_{X'}^{1/2}),\label{ccs.25}\end{equation}
 where the right hand side is defined using the induced connection
on $K_{X'}^{1/2}.$
\end{thm}
\par\noindent \bf Proof. \rm
In the proof we will identify $TX:=(TX,J)$ with 
$T^{1,0}(X,J)$ as
complex vector bundles, so that $\Lambda^{r,0}(T^{*}X,J)$ is identified
with $\Lambda^{r}(T^{*}X)$ (in particular $K_{X}=\Lambda^{n}(T^{*}X)$
and similarly for $TX':=(TX,J').$ Let us first prove (\ref{ccs.24}).
Observe that for $k$ sufficiently large, the left hand side of (\ref{ccs.24})
is given by \begin{equation}
\dim\mathcal{H}^{0,q}(X,L^{k})=(-1)^{q}\int_{X}Td(TX)\wedge 
e^{kc^{1}(L)},\label{eq:rr for hq}\end{equation}
where $Td(TX)$ is the Todd class of the complex vector bundle $(TX,J).$
Indeed, for any line bundle $L$ the Riemann-Roch theorem 
\cite{gr},\cite{gi}
applied to the complex $(\mathcal{E}^{0,*}(X,L^{k}),\overline{\partial})$
gives that the alternating sum of the dimensions of the spaces 
$\mathcal{H}^{0,j}(X,L^{k})$
is given by the right hand side in (\ref{eq:rr for hq}). Moreover,
if $L$ has index $q,$ then the dimensions of all 
$\mathcal{H}^{0,j}(X,L^{k})$
such that $j\neq q$ vanish for $k$ sufficiently large (as follows
from Proposition \ref{db1}) giving (\ref{eq:rr for hq}). Similarly,
it was shown in \cite{BoGu} that the right hand side of (\ref{ccs.24})
is given by \begin{equation}
\dim H^{0}(X',L^{k}\otimes K_{X'}^{-})=\int_{X'}Td(TX')\wedge 
e^{kc^{1}(L)+c^{1}(K_{X'}^{-})},\label{eq:rr for h0}\end{equation}
 now using the Todd class of the complex vector bundle $(TX,J').$
Using (\ref{eq:rr for hq}) and (\ref{eq:rr for h0}) and the fact that
$[X']=(-1)^{q}[X]$ as integration currents (since the orientation
depends on the almost complex structure) it is enough to show that
\begin{equation}
Td(TX)\wedge e^{kc^{1}(L)}=Td(TX')\wedge 
e^{kc^{1}(L)+c^{1}(K_{X'}^{-})}.\label{ccs.28}\end{equation}
to prove the theorem. To this end, we first recall the following basic
properties of the Todd class. Let $F$ be a complex line bundle and
$E_{1}$ and $E_{2}$ complex vector bundles over a real manifold
$X.$ Then \begin{equation}
\begin{array}{lrcl}
(i) & Td(F) & = & c^{1}(F)/(1-e^{-c^{1}(F)})\\
(ii) & Td(E_{1}\oplus E_{2}) & = & Td(E_{1})\wedge 
Td(E_{2})\end{array},\label{ccs.29}\end{equation}
where the expression in $(i)$ is to be interpreted as a formal power
series in $c^{1}(F)$, yielding a polynomial in $c^{1}(F)$, since
$c^{1}(F)^{j}$ vanishes if $j>n.$ In fact, by the {}``splitting
principle'' the properties (\ref{ccs.29}) determine $Td$
uniquely \cite{bo-tu}. Next, we will show that the following universal
identity holds \begin{equation}
Td(\overline{E})\wedge e^{c^{1}(E)}=Td(E).\label{eq:id for 
td}\end{equation}
To prove a universal identity between characteristic classes it is,
by the {}``splitting principle'' enough to prove it when $E$ is
a direct sum of line bundles over a manifold $Y.$ Moreover, by 
(\ref{ccs.29})
$(ii)$ and the multiplicativity of $e^{c^{1}}$ we may then assume
that $E$ is a line bundle. By (\ref{ccs.29}) $(i)$ the identity
(\ref{eq:id for td}) is then equivalent to the function identity \[
\frac{x}{1-e^{-x}}=\frac{-x}{1-e^{-(-x)}}\cdot e^{x}\]
which clearly holds. Let us now finish the proof of the identity 
(\ref{ccs.28}).
By the definition of $J'$ the splitting (\ref{ccs.3})
gives \[
TX'=\overline{TX_{-}}\oplus TX_{+}\]
 as complex vector bundles. Substituting this into the right hand
side of (\ref{ccs.28}) and using the multiplicative property
(\ref{ccs.29}) $(ii)$ we see that it is enough to show that
\[
Td((TX)_{-})=Td(\overline{(TX)_{-}})\wedge e^{c^{1}(K_{X'}^{-})}.\]
 Finally, since 
$c^{1}(K_{X'}^{-}):=c^{1}(\Lambda^{q}(\overline{(T^{*}X})_{-})=c^{1}((TX)_{-}),$
the identity (\ref{ccs.28}) follows from the identity (\ref{eq:id 
for td})
applied to $E=(TX)_{-}.$ This finishes the proof of (\ref{ccs.24}).
To prove (\ref{ccs.25}), note that the previous
argument also shows that (\ref{ccs.24}) remains true after
replacing $L^{k}$ by $L^{k}\otimes F$ in both sides of (\ref{ccs.24}),
where $F$ is a complex vector bundle. In particular, letting 
$F=K_{X}^{1/2}$
we get, using the decomposition (\ref{ccs.12}), that $F\otimes 
K_{X'}^{-}$
is given by\[
((K_{X}^{-})^{1/2}\otimes(K_{X}^{+})^{1/2})\otimes(K_{X}^{-})^{-1}=
K_{X^{-}}^{-1/2}\otimes K_{X^{+}}^{1/2}=K_{X'}^{1/2}\]
where we have used that $\overline{E}\simeq E^{*}:=E^{-1}$ for any
complex line bundle $E.$ This proves (\ref{ccs.25}).
\hfill{$\#$}\medskip
\begin{rem}
\rm To prove the second part of the previous theorem one could also use
that any almost complex structure whose canonical line bundle has
a square root determines a spin structure on $X$. Then the use of the
Riemann-Roch theorem may be replaced by the index theorem for the
correponding Dirac operator. In this context it is well-known that
the index only depends on the induced orientation of the real manifold
$X$. See \cite{gi}.
\end{rem}

\section{Examples: Flag manifolds}\label{sec:ex}

In this section we will recall (without giving proofs) the construction
of flag manifolds and their homogeneous line bundles, emphasizing
the complex analytical aspects. It turns out that the new almost complex
structures $J'$ (defined by (\ref{ccs.4})) in this context
are actually integrable and we show that Theorem \ref{thm:serre}
corresponds to a weak version of the Borel-Weil-Bott theorem. For
general references on flag manifolds see 
\cite{b-h}\cite{ki}\cite{fu}\cite{b-e}.
See also \cite{ku} and \cite{ka} where they are also studied from
an asymptotic point of view.

Let $K$ be a compact semi-simple real Lie group and take a maximal
connected Abelian subgroup $T$ of $K$ (i.e. a \emph{maximal torus}
of $K).$ The $K-$homogenous manifold $X:=K/T$ is called a flag
manifold. Recall that the complexification of the Lie algebra 
$\mathfrak{k}$
of $K$ decomposes as 
\begin{equation}
\mathfrak{k}_\C=\mathfrak{t}_\C\bigoplus_{\alpha\in\Delta}E_{\alpha}
\label{eq:decomp}\end{equation}
diagonalizing the adjoint action of $\mathfrak{t}$ on $\mathfrak{k}$
(acting by the Lie bracket). The label $\alpha$ of the eigen space
$E_{\alpha}$ is called a \emph{root} and it defines a non-zero element
of $\mathfrak{t}^{*}_\C :$
\[
[t,\overline{Z}_{\alpha}]=(\alpha,t)\overline{Z}_{\alpha}\]
for any element $\overline{Z}_{\alpha},$ called a \emph{root vector}, of the
\emph{root space} $E_{\alpha}.$ From (\ref{eq:decomp}) and a consistent
choice of positive roots $\Delta_{+}$ one gets a decomposition  at
the identity element $e$ of $K:$ 
\begin{equation}
T_{e}X\otimes\C\cong\Big(\bigoplus_{\alpha\in\Delta_{+}\bigsqcup-
\Delta_{+}}E_{\alpha}\Big) =:T_{e}^{1,0}X\oplus 
T_{e}^{0,1}X\label{eq:tangentbd}\end{equation}
inducing an invariant integrable complex structure on $X.$ In fact,
exponentiating the $(1,0)$ part of (\ref{eq:tangentbd}) expresses
$X$ as a holomorphic quotient, \begin{equation}
X:=K/T\backsimeq G/B,\label{Flag.3}\end{equation}
where $B$ is a \emph{Borel group} in the complexification $G$ of
$K.$ We fix a Hermitian invariant metric on $K/T$ making the
decomposition (\ref{eq:tangentbd}) orthogonal.

\par The Hermitian holomorphic line bundles on $X$ may be identified
with the \emph{weight lattice} in $\mathfrak{t}^{*}$, i.e. the elements
$\lambda$ of $\mathfrak{t}^{*}$ that exponentiate to characters
on the torus $T$, with values in $U(1)$. 
To see this, recall that in general a hermitian line
bundle over a manifold $X$ can be considered as the vector bundle
associated to a principal $U(1)-$bundle over $X.$ In our situation
$K$ is a principal $T-$bundle over $X(=K/T)$ and since $\lambda$
induces a homomorphism of the fiber $T$ into $U(1)$ it determines
a Hermitian line bundle $L_{\lambda}$ over $X.$
 The curvature two
form $\Theta_{\lambda}$ of $L_{\lambda}$ is determined by  \begin{equation}
\Theta_{\lambda}(\overline{Z}_{\alpha},Z_\beta )=\delta_{\alpha\beta}c_{\alpha}\left\langle 
\lambda,\alpha\right\rangle \label{Flag.4}\end{equation}
using the Killing form $\left\langle \cdot,\cdot\right\rangle $ on
$\mathbf{\mathfrak{t}}_{\C}^{*},$ where $\overline{Z}_{\alpha}$ is a
normalized root vector in $E_{\alpha}$ and $c_{\alpha}$ is a certain
positive number.  Formula (\ref{Flag.4}) shows that the index of
$\Theta_{\lambda}$ is equal to the {\it index} of $\lambda,$ where the
latter is defined as the number of positive roots $\alpha$ such that
$\left\langle \lambda,\alpha\right\rangle <0.$ The hyper planes
$\ker\alpha$ divide $\mathbf{\mathfrak{t}}^{*}$ into so called \emph{Weyl
chambers} and the index is constant for all $\lambda$ in the interior of a
chamber. {\it In the following we will assume that the curvature of
$L_{\lambda}$ is non-degenerate i.e. that $\lambda$ is in the interior of a
chamber.}

\begin{example}
\rm Let $K=SU(n+1).$ Then $T=U(1)^{n}$ and $G=SL(n+1,\C)$ with $B$
the subgroup of upper triangular matrices and $X$ is the manifold
of all complete flags in $\C^{n+1},$ i.e. the set of all $n-$tuples
of linear subspaces $(V_{1},...V_{n})$ such that $V_{i}$$\subsetneq 
V_{i+1}$
For example, if $n=1,$ then $X={\bf  P}^{1}$ and the conjugate complex
manifold $\overline{{\bf  P}^{1}}$ is obtained by letting $B$ be defined
by lower triangular matrices. Moreover, $\mathbf{\mathfrak{t}}=i\R$,
the weight lattice is, under proper normalization, $i\Z$ and the Weyl chambers are the positive
and negative half-axes. The element $im$ corresponds to the line
bundle $\mathcal{O}(m),$ whose sections are the homogoneous polynomials
of degree $m.$ If $n=2,$ then $X$ may be identfied with the three
dimensional manifold\begin{equation}
(Z,W)\in{\bf  P}^{2}\times{\bf  P}^{2}\ :\, 
Z_{0}W_{0}+Z_{1}W_{1}+Z_{2}W_{2}=0,\label{eq:q flag}\end{equation}
 in terms of homogenous coordinates and the action of $SU(3)$ is
given by the action \[
(A;(Z,W))\mapsto(AZ,(A^{t})^{-1}W).\]
 The weight lattice is now $i\Z^{2}$ and there are six Weyl chambers.
This follows from the representation theory of $SU(3)$ but using
the realization (\ref{eq:q flag}) it is straight forward to see that
all line bundles on $X$ are obtained as $\ 
\pi_{1}^{*}(\mathcal{O}(m))\otimes\pi_{2}^{*}(\mathcal{O}(n))$
in terms of the projections on the factors in (\ref{eq:q flag}). Moreover,
by homogenity it is, using the fiber metric induced by the Fubini-Study
metric, enough to calculate the index at a given point. Then one sees
that there are six chambers determined by linear conditions on $m$
and $n.$
\end{example}

\subsection{\label{ccswg}Change of complex structure - The 
Weyl
group}

The \emph{Weyl group} is the group generated by the reflections in
the hyper planes $\ker\alpha$ determined by the roots. It preserves
the weight lattice and acts transitively and simply on the set of
Weyl chambers. In particular, if $\lambda$ has index $q$ there is
an element $w$ of the Weyl group such that $w(\lambda)$
is positive. Dualy, the action of the Weyl group may be interpreted
as a change of the complex structure $J.$ Indeed, since \begin{equation}
\left\langle w(\lambda),\alpha\right\rangle >0\Leftrightarrow\left\langle 
\lambda,w^{-1}(\alpha)\right\rangle >0\label{Flag.6}\end{equation}
 the weight $w(\lambda)$ is in the positive Weyl chamber if and only
if the line bundle $L_{\lambda}$ is positive with respect to $J_w,$
where $J_w$ is the complex structure determined by the positive
roots $w^{-1}(\alpha)$. Hence, $L_{\lambda}$ determines a unique
invariant complex structure on $X,$ making $L_{\lambda}$ positive.
More concretely, assume that the positive roots $\alpha_{i}$ are
labeled so that \[
\left\langle \lambda,\alpha_{i}\right\rangle <0,\, i\leq 
q,\,\,\,\left\langle \lambda,\alpha_{i}\right\rangle >0,\, i>q\]
This means that the functional defined by $\lambda$ is positive precisely
on the subset \begin{equation}
\left\{ -\alpha_{1},...,-\alpha_{q},\alpha_{q+1},...\right\} 
\label{Flag.7}\end{equation}
 of the set the roots. By (\ref{Flag.6}) this set must
then be the image of the positive roots under $w^{-1}$ (which is
known to permute the roots). Furthermore, the action of $w$ induces
an isomorphism of holomorphic line bundles:\begin{equation}
\begin{array}{ccc}
L_{\mu} & \rightarrow & L_{w(\mu)}\\
\downarrow &  & \downarrow\\
G/B_{w} & \rightarrow & G/B\end{array},\label{eq:isomor}\end{equation}
where $G/B_{w}$ is the holomorphic quotient corresponding to $(X,J_w)$. 
The point is that $w$ can be identified with an element of $K$,
acting on $G$ by the adjoint action.

\subsection{The Borel-Weil-Bott theorem}

Theorem \ref{thm:serre} applied to a homogenous line bundle $L$
over the homogenous complex manifold $X$ (that can be represented
as in (\ref{Flag.3})) gives, with 
$\rho:=\frac{1}{2}\sum_{\alpha\in\Delta_{+}}\alpha:$

\begin{cor}
Assume that the weight $\lambda$ is in the interior of
a Weyl chamber and that it has index $q.$ Then, after replacing
$\lambda$ by a sufficiently large multiple, \begin{equation}
\dim H^{q}(G/B,L_{\lambda})=\dim H^{0}(G/B_{w},L_{\lambda+\rho-w^{-1}(\rho)}).
\label{Flag.9}\end{equation}
Equivalently, fixing the complex structure $J$ on $K/T:$\begin{equation}
\dim H^{q}(K/T,L_{\lambda})=\dim H^{0}(K/T,L_{w(\lambda+\rho)-\rho})\label{Flag.10}
\end{equation}
 \end{cor}
\par\noindent \bf Proof. \rm
Assume that $L_{\lambda}$ has index $q$ and let $J'$ on be the
new invariant almost complex structure determined by (\ref{ccs.4}).
Then $J'$ is an almost complex structure such that $L_{\lambda}$
is positive with respect to $J'$ and so is the complex structure
$J_w$ determined by $w$ in the Weyl group as explained in Section
\ref{ccswg}. By the uniqueness in Remark \ref{ccs3}
we have $J'=J_w.$ Hence, Theorem \ref{thm:serre}
gives (\ref{Flag.9}), but with the line bundle $L_{\lambda}\otimes 
K_{X'}^{-}$
in the right hand side. To see that $L_{\lambda}\otimes 
K_{X'}^{-}=L_{\lambda+\rho-w^{-1}(\rho)}$ note
that, given the ordering of the positive roots in Section 
\ref{ccswg},\begin{equation}
K_{X'}^{-}\leftrightarrow\sum_{i=1}^{q}\alpha_{i}=\rho-w^{-1}(\rho)\label{Flag.11}
\end{equation}
Indeed, from the definition (\ref{eq:tangentbd}) of the complex structure $J$ on $K/T$
\[
T^{1,0}X=\bigoplus_{\alpha\in\Delta_{+}}L{}_{\alpha},\]
where $L_{\alpha}$ is
the line bundle corresponding to the root $\alpha$, 
giving \[
K_{X}=\Lambda ^{n}(T^{*1,0}X)\simeq\bigotimes_{\alpha\in\Delta_{+}}L{}_{-\alpha}=L_{-2\rho}\]
and a similar argument gives the first correspondence in (\ref{Flag.11}).
Finally, since the image of the positive roots under $w^{-1}$ is
given by (\ref{Flag.7}), \[
\rho-w^{-1}(\rho)=\frac{1}{2}(\sum_{i=1}^{q}\alpha_{i}+\sum_{i=q+1}^{n}\alpha_{i})-\frac{1}{2}(\sum_{i=1}^{q}-\alpha_{i}+\sum_{i=q+1}^{n}\alpha_{i})=\sum_{i=1}^{q}\alpha_{i}\]

Now the induced isomorphism (\ref{eq:isomor}) applied to 
$\mu=\lambda+\rho-w^{-1}(\rho)$
proves (\ref{Flag.10}).
\hfill{$\#$}\medskip

\par 
The previous corollary (in the formulation (\ref{Flag.10})) is
a weak version of Bott's generalization of the Borel-Weil theorem
\cite{bo},\cite{bo2}. The Borel-Weil-Bott theorem may also be proved
using Lie algebra cohomology \cite{ko}\cite{z}.

\section{\label{ax}Appendix: The affine bundle $AX$}

We will define an affine bundle $AX$ over $X$ with symplectic form
$\Omega$ so that the global sections of $AX$ are the unitary connections
of the Hermitian line bundle $L$ over $X.$ Given an open set $U$
and a local unitary frame $t$ for $L$ over $U$ we identify $(AU,\Omega)$
with $(T^{*}U,dp\wedge dx)$ in terms of the usual coordinates $(x,p)$
on $T^{*}U.$ If $\widehat{t}=e^{ig}t$ is another unitary section
the two identifications are assumed to be related by\[
(x,\widehat{p}_{1},...,\widehat{p}_{n})=(x,p_{1}-\frac{\partial}{\partial 
x_{1}}g,...,p_{n}-\frac{\partial}{\partial x_{n}}g),\]
 Hence, $\Omega=dx\wedge dp$ is a globally well-defined symplectic
two-form on $AX.$ Given a global connection $\nabla$ represented
by $d+A$ with respect to the frame $t$ the transformation property
(\ref{ccs.11}) now shows that $(x,iA_{1}(x),..,iA_{n}(x))$
defines a global section of $AX.$ 

Notice that the local characteristic variety $\Sigma$ in Proposition
\ref{db2} corresponds globally to the graph in $AX$ of the canonical
connection $\nabla$ on the Hermitian holomorphic line bundle $L.$
Indeed, by \no{db.22.5} and \no{db.34} we get locally on $\Sigma$\[
pdx=\textrm{Re$(\frac{2}{i}\partial\phi)$}=i(-\partial\phi+\overline{\partial}\phi)\]

By (\ref{ccs.15}) the right hand side equals $iA,$ where
$A$ is the local one form associated to $\nabla$ with respect to
$t=e^{\phi}s.$

Finally, for comparison with \cite{g-u}\cite{l} observe that any
given unitary connection $\nabla$ on $L$ induces a global 
isomorphism\begin{equation}
\Phi_{\nabla}:\, AX\leftrightarrow 
T^{*}X,\,\,\,(x,p)\mapsto(x,p_{1}-iA_{1},...),\label{ax.1}\end{equation}
The map is defined using local frames $t$ as above and it maps the
graph of the section of $AX$ corresponding to $\nabla$ to the zero-section
in $T^{*}X.$ We get that 
$$(\Phi_{\nabla}^{-1})^{*}(\Omega)=d(-\gamma)+\pi^{*}(-i\Theta),$$
where $\gamma$ is the tautological $1-$form on $T^{*}X$ and 
$\pi^{*}(-i\Theta)$
is the normalized curvature of $\nabla$ pulled back from $X.$ The
bundle $AX$ may also be defined by symplectic reduction of 
$(T^{*}Y,d(-\gamma)$
where $Y$ is the unit circle bundle in $L^{*}$ (compare the proof
of Theorem 2.3 in \cite{s-z}).

\end{document}